\documentclass{amsart}
\usepackage{graphics}
\usepackage{rlepsf}
\usepackage{a4}
\usepackage{amssymb}

\newcommand{\C}{\mathbb C}
\newcommand{\Z}{\mathbb Z}
\newcommand{\cL}{\mathcal L}
\newcommand{\cG}{\mathcal G}\newcommand{\LG}{\Lambda \cG}
\newcommand{\cH}{\mathcal H}
\newcommand{\oG}{\overline G}
\newcommand{\cB}{\mathcal B}
\newcommand{\cE}{\mathcal E}

\newcommand{\GG}{\mathfrak G}
\newcommand{\id}{\text{Id}}
\newcommand{\Gconj}{G^{\text{c}}}

\newcommand{\Vect}{\mathcal{V}\text{ect}}

\newcommand{\vpic}[1]{}

\newcommand{\nattrans}{\Rightarrow}
\newcommand{\pretheta}{{}^\theta\!}

 \DeclareMathOperator{\Rep}{Rep}
 \DeclareMathOperator{\thetaRep}{\Rep^\theta\!}
 \DeclareMathOperator{\cRep}{\mathcal Rep}
 \newcommand{\thetacRep}{\cRep^\theta\!}
 \DeclareMathOperator{\Fun}{\mathcal Fun}
 \DeclareMathOperator{\Hom}{Hom}
 \DeclareMathOperator{\Ob}{Ob}
 
 \DeclareMathOperator{\Tr}{Tr} 
 \DeclareMathOperator{\End}{End} \DeclareMathOperator{\Image}{Im}
 \DeclareMathOperator{\Maps}{Maps} \DeclareMathOperator{\Mor}{Mor}
 \DeclareMathOperator{\Aut}{Aut} \DeclareMathOperator{\Par}{Par}
 \DeclareMathOperator{\Groth}{Groth}
 \DeclareMathOperator{\rank}{rank}

\newtheorem{thm}{Theorem}
\newtheorem{prop}[thm]{Proposition}
\newtheorem{cor}[thm]{Corollary}
\newtheorem{lemma}[thm]{Lemma}

\begin{document}

\title[The twisted Drinfeld double via gerbes and groupoids]
      {The twisted Drinfeld double of a finite group\\
       via gerbes and finite groupoids}
\author{Simon Willerton}
\address{University of Sheffield}
\email{S.Willerton@sheffield.ac.uk}

\begin{abstract}
The twisted Drinfeld double (or quasi-quantum double) of a finite
group with a 3-cocycle is identified with a certain twisted
groupoid algebra.  The groupoid is the loop (or inertia) groupoid
of the original group and the twisting is shown geometrically to
be the loop transgression of the 3-cocycle.  The twisted
representation theory of finite groupoids is developed and used to
derive properties of the Drinfeld double, such as representations
being classified by their characters.

This is all motivated by gerbes and 3-dimensional quantum field
theory.  In particular the representation category of the twisted
Drinfeld double is viewed as the `space of sections' associated to
a transgressed gerbe over the loop groupoid.
\end{abstract}
\maketitle %
\thispagestyle{empty}
\section*{Introduction}

There were several motivations for this work.  The first was to
understand what the twisted Drinfeld double of a group
\emph{really} was, and in particular to understand exactly what
the twisting was.  The second was to see how well Dijkgraaf-Witten
quantum field theory worked in the two-gerbe paradigm.  The third
was to understand the Freed-Hopkins-Teleman result on twisted
equivariant K-theory in this baby example.  Actually what I really
want to understand is Chern-Simons theory and the hope was this
might some give some insight in that direction.

The twisted Drinfeld double of a finite group arose in the work of
Dijkgraaf, Pasquier and Roche in relation to orbifold conformal
field theory and was shown by Freed
\cite{Freed:HigherAlgebraicStructures} (following Freed and Quinn
\cite{FreedQuinn} and Altschuler and Coste
\cite{AltschulerCoste:QuasiQuantum} to be the relevant quantum
group for the quantum field theory from a finite group twisted by
the cocycle. Given a finite group $G$ and a three-cocycle
$\omega\in Z_{\text{gp}}^3(G,U(1))$ the twisted Drinfeld double
$D^\omega(G)$ can be defined to be the algebra over the complex
numbers, with a basis $\{\langle \stackrel{g}\leftarrow x
\rangle\;|\; g,x\in G\}$ and with the product defined by
 \[\langle \stackrel h \leftarrow y\rangle \cdot
  \langle \stackrel g \leftarrow x\rangle
   :=
   \delta_{y,gxg^{-1}}
      \frac{\omega(h,g,x)\omega(hgx(hg)^{-1},h,g)}{\omega(h,gxg^{-1},g)}
       \langle \stackrel{hg}\longleftarrow x \rangle.
\]
Freed showed that in the untwisted case this can be thought of as
the \emph{groupoid algebra} of a certain finite groupoid.  In this
paper I develop that idea considerably further and show that the
scalar factor in the definition of the product --- the twisting
--- is essentially a two-cocycle $\tau(\omega)$ on the \emph{loop groupoid}
$\Lambda\oG$, and that the twisted Drinfeld double is the
corresponding twisted groupoid algebra.

Following that, I then show that many properties of the twisted
Drinfeld double
--- such as its semi-simplicity, its character theory and the
decomposition of its representation group --- follow from the fact
that it is a twisted groupoid algebra.  (Some of these properties
have been algebraically presented in the somewhat terse
\cite{AltschulerCosteMaillard}.)  The twisted representation
theory of finite groupoids, which is developed here, is hardly
more complicated than the twisted (or projective) representation
theory of finite groups which, in turn, is hardly more complicated
than the ordinary representation theory of finite groups.  In fact
looking at the twisted representation theory of finite groupoids
throws some more light on the twisted representation theory of
finite groups, so that the notion of twisted characters is just an
example of the transgression phenomenon.

From the point of view of gerbes the key topological construction
is the Parmesan map.  I show how the $S^1$--transgression map is
expressed in general for $n$--gerbes on finite groupoids.  To be
more specific, by an $n$--gerbe on a finite groupoid I am really
meaning a $U(1)$--valued $(n+1)$--cocycle on $\cG$.  Associated to
$\cG$ is the \emph{loop} groupoid (aka.\ the \emph{inertia}
groupoid) and the transgression map associates an $n$--cocycle on
the loop groupoid. This is the Parmesan construction of Section~1
and it is that which allows us to move up and down the gerbe
ladder.

From the point of view of quantum field theory we see how
Dijkgraaf-Witten theory fits into the two-gerbe philosophy to
explain this it is useful to have a digression on $n$--gerbes and
transgression.

\subsection*{TQFTs and $n$--gerbes}
Although not strictly necessary for understanding the results of
this paper the following will give some idea of the thinking
behind them.  

\subsubsection*{Topological quantum field theory}
Suppose $G$ is a compact group.  If $G$ is finite then we will be
considering Dijkgraaf-Witten theory, if $G$ is simply connected
and simple then we will be considering Chern-Simons theory.  Given
a ``level'' $k\in H^4(BG,\Z)$ we can construct a three-dimensional
topological quantum field theory. This means, vaguely, to each
closed three-manifold we associate a complex number, to each
closed two-manifold we associate a vector space and to each three
manifold with boundary we associate a vector in the vector space
associated to its boundary.  In fact we might want to go further
and to each closed one-manifold associate a category and even
associate to each closed zero-manifold a two-category. The
physical way to do this (at least at the vector space and vector
level) is using path integrals; so for instance the invariant of a
closed three-manifold $M$ is the integral of a certain action over
the space of connections on a $G$ bundle over $M$, and the vector
space associated to a two-manifold, $\Sigma$ is a certain space of
sections of a line bundle over the ``moduli space'' of flat
$G$--connections on $\Sigma$.  The philosophy of the current
paper, which I learnt from Freed (see for example
\cite{Freed:HigherAlgebraicStructures}) and Brylinski and
McLaughlin \cite{BrylinskiMcCaughlin}, is that the hierarchy of
$n$-gerbes gives a mathematical framework for trying to understand
these path integrals.

\subsubsection*{Transgression and $n$--gerbes}
Implicit in Dijkgraaf and Witten's work is that it is sensible to
interpret the level $k\in H^4(BG,\Z)$ as a two-gerbe with
connection.  Being slightly vague, a two-gerbe is a geometric
representative of a class in $H^4(X,\Z)$, just as a principal
$U(1)$--bundle is a geometric representative of a class in
$H^2(X,\Z)$.  There are many ways of representing a class in a
`geometric' way, but for the sake of this discussion, we can take
a two-gerbe to be a sort of bundle of two-categories. There is a
hierarchy of $n$--gerbes, which we can indicate as follows:

\begin{center}
\begin{tabular}{c|c|c} $n$&``$n$--gerbe''& represents a class
in\dots\\ \hline
$-2$&$\Z$--valued function &$H^0(X,\Z)$\\
$-1$&$U(1)$--valued function &$H^1(X,\Z)$\\
$0$&principal $U(1)$--bundle&$H^2(X,\Z)$\\
$1$&``bundle of categories''&$H^3(X,\Z)$\\
$2$&``bundle of two-categories''&$H^4(X,\Z)$
\end{tabular}
\end{center}
In the geometric ways of representing higher gerbes there are
appropriate notions of connections, generalizing that on
$U(1)$--bundles, and an $n$--gerbe with connection, essentially by
definition, represents a class in the more subtle Deligne
cohomology group.  For $n$--gerbes with connection there are
notions of holonomy and transgression, the basic example of which
is got from a $U(1)$--bundle-with-connection on a space $X$.  The
holonomy of this can be thought of as a $U(1)$--valued function on
the loop space of $X$, so there is a way of going from a $0$-gerbe
(with connection) on $X$ to a $(-1)$--gerbe (with connection) on
$\Maps (S^1,X)$. Similarly given a one-gerbe on $X$, this can be
transgressed to a $U(1)$--bundle with connection on the loop space
or to a $U(1)$--valued function on the space of maps of a surface
$\Sigma$ into $X$.  More generally, for an oriented, closed
$m$--manifold $Y$, with $0\le m \le n+2$ there is a transgression
map
\begin{multline*}
  \{\text{$n$--gerbes with connection on $X$}\}\\
  \longrightarrow
  \{\text{$(n-m)$--gerbes with connection on $\Maps(Y,X)$}\},
\end{multline*}
which on the level of cohomology is just the composition
\[H^{n+2}(X,\Z)\stackrel{\text{ev}^*}\longrightarrow
H^{n+2}(\Maps(Y,X)\times Y,\Z)\stackrel{\int_Y}{\longrightarrow}
H^{(n-m)+2}(\Maps(Y,X),\Z),
\]
where $\text{ev}\colon \Maps (Y,X)\times Y\to X$ is just the
evaluation map.  This is constructed on the level of Deligne
cohomology by Gomi and Terashima
\cite{GomiTeramshima:FiberIntegration}.

It its possible to make sense of transgression for $m$--manifolds
with boundary and that is essential in setting up the full quantum
field theory but it is not something that we will be using here.
The next thing to do is to consider the `spaces of sections'.

\subsubsection*{Spaces of sections of gerbes}
In order to build a topological quantum field theory we need to
understand the ``space of sections'' of an $n$--gerbe, which is a
higher analogue of geometric quantization.  The term ``space of
sections'' is not entirely adequate, and perhaps ``integral'' or
``push-forward to a point'' would be better.  First consider the
case of a $U(1)$--bundle, to obtain its space of sections we must
tensor with the complex line to obtain a line bundle and then we
take the sections of that, to obtain a vector space.  What we
really want at this point is a finite-dimensional vector space and
this is where we get to the nub of geometric quantization, the
thorny issue of polarization. Finding a polarization results in
cutting down drastically the number of sections that we consider,
for instance if a holomorphic structure exists on the base and the
bundle, then we can restrict to holomorphic sections, obtaining a
finite dimensional space of section in the case of a compact base.
Ignoring the issue of polarization for a moment, the correct
notion of space of sections for a $U(1)$--valued function is just
the integral of the function, so that is a complex number.  In the
case of a gerbe, the ``space of sections'' (after the analogue of
tensoring with the complex numbers) forms a linear category.  The
process is polarization for gerbes seems even more mysterious.
Fortunately polarization is not an issue that has to be addressed
in this paper as the base spaces are all discrete.  So we can make
sense of the space of sections of a gerbe.

Modulo the difficulties of defining the spaces of polarized
sections, given an $n$--gerbe over a space $X$ we thus have a
philosophy of obtaining a linear $(n-m)$--category from a closed,
oriented $m$--manifold $Y$, namely take the space of (polarized)
sections of the complexified gadget associated to the transgressed
$(n-m)$--gerbe over $\Maps(Y,X)$.

\subsubsection*{Stacks, groupoids and gerbes}
I actually want to work with `stacks' rather than spaces.  This
means that I want to take extra symmetry information into account
and means in practice dealing with groupoids rather than spaces.
So the points in the space are now objects of a groupoid and there
is the extra information of the morphism sets between points.  The
example here is the stack of $G$--bundles on a compact manifold
$M$, where $G$ is a finite group.  (As $G$ is finite, these will
necessarily be flat bundles.)  To get a groupoid representing this
stack we can take a finite-pointed model $\Pi_1(M)$ for the
fundamental groupoid, by which I mean a groupoid whose objects
correspond to a finite number of points on $M$, with at least one
point on each connected component of $M$ and whose morphisms are
the homotopy classes of paths between these points; in particular
if $M$ is connected and $m_0$ is a base-point then we could take
the groupoid with one object and with the fundamental group
$\pi_1(M,m_0)$ as the group of morphisms.  It is convenient to
allow more than one point on each component, as this makes things
more natural if we are considering spaces with boundaries.  Now a
groupoid representing the stack of flat $G$--bundles on $M$ is
given by the groupoid of functors $\GG(M):=\Fun(\Pi_1(M),\oG)$,
where $\oG$ is the groupoid with one object and $G$ as its group
of morphisms, while the functor groupoid has functors as its
object and natural transformations as its morphisms.  Note two
things here, firstly $\GG(M)$ is a finite groupoid as $\Pi_1(M)$
is finitely generated and $\oG$ is finite. Secondly, $\Pi_1(M)$ is
only defined up to equivalence, and hence so is $\GG(M)$.

To be concrete, here an $n$--gerbe on a finite groupoid will mean
a $U(1)$--valued, $(n+1)$--cocycle on the groupoid.  We can define
transgression similar to above, so that in particular if $\omega$
is a $3$--cocycle on $G$ then for $M$ a closed (triangulated)
$m$--manifold there is the transgression $\tau_M\omega$ which is a
$(2-m)$--gerbe on $\GG(M)$.  (The cocycle will actually depend on
the triangulation of $M$, but it's cohomology class will not.)  We
can associate a ``complexified bundle'' $(\tau_M\omega)_\C$ and
take its space of sections to get a $(2-m)$--category
$V(M):=\Gamma_{\GG(M)}((\tau_M\omega)_\C)$ which is the invariant
of $M$.

Thus given $G$ and $\omega\in Z^3_{\text{gp}}(G,U(1))$ we get the
sequence of assignments
 \[\begin{array}{ccccccc}\parbox{4em}{\center\text{manifold:}\\[.3em]
    \footnotesize$M^m$}
   &\leadsto&\parbox{15em}{\center\text{$(2-m)$--gerbe over finite
   groupoid:}\\[.3em]
     \footnotesize $\tau_M\omega$ over $\GG(M)=\Fun(\Pi_1(M),\oG)$}
 &\leadsto&\parbox{10em}{\center\text{space of sections:}\\[.3em]
  \footnotesize $V(M):=\Gamma_{\GG(M)}((\tau_M\omega)_\C)$}\\
 \end{array}\]

\subsubsection*{Dimension and $S^1$-transgression}
Given a linear $n$--category, you can form the Grothendieck
$(n-1)$ category of `equivalence classes', and by the
``dimension'' of an $n$-category I mean some sort of
complexification of the Grothendieck $(n-1)$--category --- I only
really know what this means for $n=1$ and $n=0$, see below.  From
the point of view of topological quantum field theory, we expect
the $(1-m)$--category $V(S^1\times M)$ to be the dimension of
$V(M)$ the $(2-m)$--category associated to $M$. But if $\Pi_1(M)$
is a finite model for the fundamental groupoid of $M$ then
$\overline \Z\times \Pi_1(M)$ is a finite model for the
fundamental groupoid of $S^1\times M$.  For any finite groupoid
$\cG$ we can define the \emph{loop groupoid} (or \emph{inertia
groupoid}) $\Lambda \cG$ to be $\Fun(\overline \Z, \cG)$. Then
  \begin{align*}\GG(S^1\times M)
    &=\Fun(\Pi_1(S^1\times M),\oG)
    =\Fun(\overline \Z\times \Pi_1(M),\oG)\\
    &=\Fun(\overline \Z,\Fun(\Pi_1(M),\oG))
    =\Lambda\GG(M).
  \end{align*}
Furthermore as
 $\tau_{(S^1\times M)}\omega=\tau_{S^1}\tau_M\omega$,
we can write
\[V(S^1\times M)
             =\Gamma_{\Lambda \GG(M)}(\tau_{S^1}(\tau_M\omega)_\C)\]
and so we expect this to be the dimension of
$V(M)=\Gamma_{\GG(M)}((\tau_M\omega)_\C)$.

Actually, given what we will see below, we might expect something
more general to be the case.  Thinking of our model of gerbes to
be simplicial cocycles, in Section~? we define the general
$S^1$--transgression map  $\tau\colon Z^n(\cG,U(1))\to
Z^{n-1}(\Lambda \cG,U(1))$, then we will see in the sections below
that for any finite groupoid $\cG$ and any one- or two-cocycle
$\xi$ on $\cG$ we have
 \[
   \dim(\Gamma_{\cG}(\xi_\C))=\Gamma_{\Lambda \cG}(\tau(\xi)_\C).
   \eqno{(\star)}
 \]
At this stage this just looks like a not-very meaningful
collection of symbols, so it is worth unpacking it a little in the
case that $\xi$ is of degree one or two.

When $\xi=\alpha$ is a one-cocycle on $\cG$ then $\alpha$ can be
interpreted as a $U(1)$--bundle on $\cG$, and the space of
sections $\Gamma_{\cG}(\alpha)_\C$ is the space of (flat) sections
of the associated complex line bundle. The transgression
$\tau(\alpha)$ is a locally constant $U(1)$--valued function on
the loop groupoid $\Lambda \cG$ and the ``space of sections''
$\Gamma_{\Lambda\cG}\tau(\alpha)_\C$ is just the integral of
$\tau(\alpha)$ over the loop groupoid, and does indeed give the
dimension of the vector space $\Gamma_{\cG}(\alpha)_\C$.  In this
case we can write ($\star$) as
 \[\dim(\text{sections of the line
   bundle $\alpha_\C$ over $\cG$})=\int_{\Lambda \cG}\tau(\xi),\]
and
this is Theorem~\ref{Theorem:DimSectionsBundle}.

When $\xi=\theta\in Z^2(\cG,U(1))$ is a two-cocycle, then the
gerbe $\theta$ can be thought of as an central extension of $\cG$
by $U(1)$, and the ``space of sections'' $\Gamma_{\cG}(\theta)_\C$
is the category of $\theta$--twisted representations of $\cG$,
whereas $\Gamma_{\Lambda \cG}(\tau(\theta)_\C)$ is the vector
space of $\theta$--twisted characters, so the assertion comes down
to equivalence classes of twisted representations being classified
by twisted characters.  In this case ($\star$) can be written as
 \[\Groth(\thetacRep(\cG))\otimes \C\cong \{\text{sections of the line
   bundle $\tau(\theta)_\C$ over $\Lambda\cG$}\},\]
and this is Theorem~\ref{Theorem:DimSectionsGerbe}.

Note that when the groupoid $\cG$ is just a group $\oG$ and the
two-cocycle $\theta$ is trivial, the above is saying that the
characters of $G$ classify the representations and span the class
functions, in other words the character map induces an isomorphism
$$\Rep(G)\otimes\C\cong
  \{\text{class functions}\ G\to\C\}.$$

In order to make sense of ($\star$) when $\xi$ is a higher order
cocycle (or gerbe) it is necessary to understood higher categories
better than I currently do.

\subsubsection*{Dijkgraaf-Witten theory and the Drinfeld double}
We can see how the above ideas of transgression and dimension fit
into the picture of Dijkgraaf-Witten theory in low dimension.  The
following table gives an indication of what is going on.  We
assume a fixed finite group $G$ and three-cocycle $\omega\in
Z_{\text{gp}}^3(G,U(1))$.

\begin{center}
\begin{tabular}{c|c|c|c}
  \parbox[c][1.5\height]{4em}{\center{manifold\\$M$}} &
  \parbox[c][1.5\height]{4em}{\center{groupoid\\$\GG(M)$}} &
  \parbox[c][1.5\height]{7em}{\center{gerby object\\$\tau_M\omega$}} &
  \parbox[c][1.5\height]{12em}{\center{``space of sections''\\
                      $\Gamma_{\GG(M)}\bigl((\tau_M\omega)_\C\bigr)$}}
  \\
\hline%
%
   $S^1$  &
   $\Lambda\oG$ &
   \parbox[c][1.5\height]{7em}{\center{$\tau(\omega)$\\gerbe}} &
   \parbox[c][1.5\height]{11em}{\center{category of representations of $D^\omega(G)$}}
 \\
   $\mathbb{T}^2$ &
   $\Lambda^2\oG$ &
   \parbox[c][1.5\height]{7em}{\center{$\tau^2(\omega)$\\$U(1)$--bundle}} &
   \parbox[c][1.5\height]{11em}{\center{vector space of $\omega$--twisted
                                           elliptic characters of $G$}}
 \\
   $\mathbb{T}^3$ &
   $\Lambda^3\oG$ &
   \parbox[c][1.5\height]{7em}{\center{$\tau^3(\omega)$\\ $U(1)$--function}} &
   \parbox[c][1.5\height]{11em}{\center{number of $\omega$--twisted elliptic characters of $G$}}
\end{tabular}
\end{center}
%
By an $\omega$--twisted elliptic character of $G$ I mean a
complex-valued function $\chi$, defined on pairs of commuting
elements of $G$, which satisfies the following conjugation
relation for any $h\in G$:
  \[
  \chi(hgh^{-1},hxh^{-1})
    =\frac{\omega(h,x,g)\omega(hgh^{-1},h,x)\omega(hxh^{-1},hgh^{-1},h)}
          {\omega(h,g,x)\omega(hxh^{-1},h,g)\omega(hgh^{-1},hxh^{-1},h)}
           \chi(g,x).
  \]

Note that in the final column of spaces of sections, by the two
instances of the dimension relation ($\star$) alluded to above,
moving from the bottom up one row and from the second to bottom up
one row are both processes of categorification.

\subsection*{What is in this paper}
Section~1 is the geometric part of the paper.  It starts at the
beginning with preliminaries on classifying spaces of groupoids
and cocycles on groupoids.  The intention is to make this
accessible to people who are not familiar with these ideas.
Section~1.3 is the heart of the geometric part, in it the loop
groupoid of a groupoid is defined, and the Parmesan map, between
the classifying space of the loop groupoid and the free loop space
on the classifying space of the original groupoid, is constructed.
Theorem~\ref{Thm:GroupoidParmesan} states that this is a homotopy
equivalence. This is then used to give an explicit formula for the
loop transgression map on the level of cocyles, $\tau\colon
Z^n(\cG,U(1))\to Z^{n-1}(\Lambda \cG,U(1))$, which can be thought
of as mapping (flat) $(n-1)$--gerbes on a groupoid to (flat)
$(n-2)$--gerbes on its loop groupoid.


Section~2 is the algebraic part of the paper. This is where the
theory of twisted representations of groupoids is developed from
the point of view of gerbes.  First of all there is a quick
reminder on the theory of twisted representations of finite
groups.  Then for a finite groupoid $\cG$, considering $n=0,1,2$
in turn, I show how a normalized $n$--cocycle in $Z^n(\cG,U(1))$
can be thought of as an $(n-1)$--gerbe on $\cG$, ie.\ respectively
a $U(1)$--valued function, a principal $U(1)$--bundle and a
$U(1)$--gerbe, and how its `complexification' has a `push-forward
to a point' which is an $(n-1)$--category, these are respectively
the integral over $\cG$, the space of sections of the associated
complex line bundle, and the category of twisted representations.
Furthermore, for $n=1,2$, I show that the decategorification of
the `space of sections', ie.\ the dimension or the
complexification of the Grothendieck group, is given by the space
of sections the transgressed $(n-2)$--gerbe over the loop
groupoid.  For $n=2$ this generalizes the fact that the space of
class functions on a finite group is the complexification of the
Grothendieck group of the category of representations..  The final
subsection shows how the group of twisted representations of a
groupoid decomposes into the twisted representation groups of its
automorphism groups of its components --- this is essentially a
version of Dijkgraaf-Pasquier-Roche induction.

Section~3 is where the above work is put to use in the context  of
the twisted Drinfeld double of a finite group.  For $G$ a finite
group and $\omega\in H^3(G,U(1))$ the `level' or `twisting', it is
shown that $D^\omega(G)$ is just ${}^{\tau(\omega)}\C \Lambda\oG$
the twisted groupoid algebra of the loop groupoid on $G$.  The
representations of this are then the same as the
$\tau(\omega)$--twisted representations of the loop groupoid
$\Lambda\oG$ and all of the previous machinery can be applied to
see such things as the semi-simplicity of $D^\omega(G)$, DPR
induction, the fact that its representation group is a twisted
equivariant K-group ${}^{\tau(\omega)}K_G(\Gconj)$, and the fact
that its representations are classified by their characters, the
characters being identifiable with `twisted elliptic characters on
$G$'.

\subsection*{Acknowledgements}
I would like to thank Olivia Rossi-Doria for her nomenclature, and
the following people for various enlightening conversations:
 Mark Brightwell, J\o rgen Ellegaard Andersen, Dan Freed, John Greenlees, Kirill
 Mackenzie, Ieke Moerdijk, Michael Murray, Elmer Rees, Justin Roberts,
 Joost Slingerland,
 Neil Strickland, and Paul Turner.  I would also like to point out
 that some similar ideas on loop groupoids were explored to different ends by
 Lupercio and Uribe \cite{LupercioUribe:GerbesTwistedKtheory}.
\section{Groupoids, classifying spaces and cohomology}
The goal of this section is to give a formula for the
loop-transgression map $Z^*(\cG,U(1))\to
Z^{*-1}(\Lambda\cG,U(1))$.  This describes the twisting in the
twisted Drinfeld double and also is useful in describing
characters of twisted representations of groups.

In Section~\ref{Section:ClassifyingSpace} the classifying space of
a category is introduced using the geometric bar construction,
rather than the usual common barycentric version, as this is more
useful in the formulas derived later.  In
Section~\ref{Section:Equivalence} an explicit formula for the
well-known fact that a groupoid is equivalent to the union of its
automorphism groups is given, as this will be useful.  In
Section~\ref{Section:LoopGroupoidsParmesan} the loop groupoid (or
inertia groupoid) of a finite groupoid and the Parmesan
construction are introduced, the latter being a homotopy
equivalence $\cB\Lambda\cG\to \cL\cB\cG$ from the classifying
space of the loop groupoid to the free loop space of the
classifying space of the original groupoid. From this it is
straight forward to write down the transgression map
$Z^n(\cG,U(1))\to Z^{n-1}(\Lambda\cG,U(1))$. To prove that the
Parmesan map is a homotopy equivalence it suffices to prove it for
finite groups, this is the content of
Section~\ref{Section:ParmesanForGroups}. The key here is that the
loop groupoid $\Lambda\oG$ is also the action groupoid of $G$ on
itself by conjugation.

\subsection{The classifying space of a groupoid}
\label{Section:ClassifyingSpace}%
Recall firstly that a finite groupoid is a category with a finite
number of objects and morphisms, in which all of the morphisms are
invertible.  A good example to bear is mind has one object, and
its set of morphisms form a finite group $G$: this groupoid will
be denoted $\oG$.

\begin{figure}
$${\relabelbox\small
\epsfbox{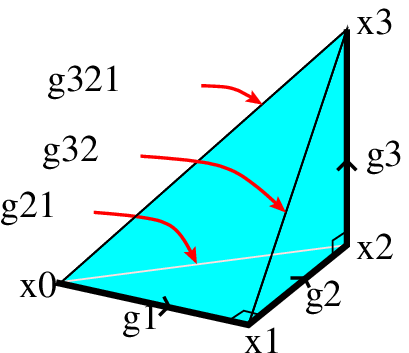} \relabel {g1}{$g_1$} \relabel{g2}{$g_2$}
\relabel {g3}{$g_3$} \relabel {x0}{$x_0$} \relabel{x1}{$x_1$}
\relabel{x2}{$x_2$} \relabel {x3}{$x_3$} \relabel{g21}{$g_2\circ
g_1$} \relabel{g32}{$g_3\circ g_2$} \relabel{g321}{$g_3\circ
g_2\circ g_1$}
\endrelabelbox}$$
\caption{The three-simplex in the classifying space $\cB\cG$
corresponding to $x_3\stackrel{g_3}\leftarrow
x_2\stackrel{g_2}\leftarrow x_1\stackrel{g_1}\leftarrow x_0$}
\label{Figure:HillTetrahedron}
\end{figure}

Given a finite groupoid $\cG$, we will associate to it a
topological space $\cB\cG$, its classifying space.   This is
formed by gluing simplices together.  (I will be thinking
geometrically and blur the distinction between a simplicial
complex and its geometric realization.)  The $n$--simplices of the
classifying space $\cB\cG$ are labelled by the strings of $n$
composable morphisms $x_n\stackrel{g_n}\leftarrow
\ldots\stackrel{g_2}\leftarrow x_1\stackrel{g_1}\leftarrow x_0$.
(Note that I sometimes write arrows ``backward'', the reason being
that notationally it makes composition behave better.)  The
corresponding $n$--simplex is thought of as an $n$--simplex in
$n$--space with the vertices labelled by $x_0, \ldots, x_n$, with
the directed edge from $x_{i-1}$ to $x_i$ of length one, parallel
to the $i$th axis and  labelled by $g_i$, with the other edges
labelled so as to make the two-faces commutative triangles. A
three-simplex is pictured in Figure~\ref{Figure:HillTetrahedron}.
The classifying space $\cB\cG$ is then obtained by gluing common
faces together.

This means that a point in the classifying space is given by a
string $x_n\stackrel{g_n}\leftarrow \ldots\stackrel{g_2}\leftarrow
x_1\stackrel{g_1}\leftarrow x_0$ and some coordinates $0\le
t_n\le\dots\le t_1\le 1$ in the simplex; this point will be
denoted
  $$[g_n,t_n|\ldots|g_2,t_2|g_1,t_1]
  .$$
The gluing of faces means that there is the following
identification of points in $\cB\cG$.
\begin{align*}
[g_n,t_n|g_{n-1},t_{n-1}|\ldots|g_1,1]&=[g_n,t_n|\ldots|g_2,t_2];\\
[g_n,t_n|\ldots|g_{i+1},t|g_{i},t|\ldots|g_1,t_1]
 &=[g_n,t_n|\ldots|g_{i+1}\circ g_i ,t|\ldots|g_1,t_1];\\
[g_n,0|g_{n-1},t_{n-1}|\ldots|g_1,t_1]
  &=[g_{n-1},t_{n-1}|\ldots|g_1,t_1].
\end{align*}
This is just saying that a point is in the face of a simplex when
two adjacent coordinates are equal, when the first coordinate is
$1$ or when the final coordinate is $0$.

I think it is more common to see barycentric coordinates used for
simplices, but the coordinate system given here will have
advantages when prisms are decomposed below.  Another common
variation is the imposition of degeneracy conditions, so that an
$n$--simplex in which one of the morphisms is an identity
morphism, is identified with the $(n-1)$--simplex obtained by
removing that morphism.  I do not impose these conditions as it is
not necessary.  However, when considering cocycles it will be
convenient to consider \emph{normalized} ones, these are the
cocycles which vanish on simplices containing an identity morphism
(see below).  The normalized and non-normalized cochain complexes
are homotopic and so have the same cohomology groups, thus it is
primarily a matter of convenience as to which is used.

 Any functor
$F\colon \cG\to\cH$ between finite groupoids determines an obvious
map $\cB F\colon \cB\cG\to \cB\cH$ between classifying spaces
given by
  $$\cB F\colon [g_n,t_n|\ldots|g_2,t_2|g_1,t_1]
    \mapsto [F(g_n),t_n|\ldots|F(g_2),t_2|F(g_1),t_1].$$

By definition, the cohomology of $\cG$ is taken to be the
simplicial cohomology of $\cB\cG$ (thought of as a simplicial
complex), so an $n$--cochain on $\cG$ means a function from the
set of $n$--simplices of $\cB\cG$, which will be normalized so
that $[g_n|\dots|\id_{x_i}|\dots|g_1]$ maps to the unit.
Furthermore $\omega$ is an $n$--cocycle if it is an $n$--cochain
which vanishes on boundaries, so that
  $$\omega\bigl([g_{n-1}|\ldots|g_1]\bigr)
    \omega\bigl([g_ng_{n-1}|\ldots|g_1]\bigr)^{-1}
    \cdots
    \omega\bigl([g_n|\ldots|g_2g_1]\bigr)^{\pm 1}
    \omega\bigl([g_n|\ldots|g_2]\bigr)^{\mp 1}
    =1.
  $$

If $G$ is a finite group then $\oG$ was defined to be the category
with one object, with the morphism set being $G$ and with
composition being multiplication in $G$.  Then $\cB\oG$ is a model
for the classifying space for $G$ so that $H^*(\oG)$ is the usual
group cohomology of $G$.  This particular model is known as the
geometric bar construction (see \cite{Gajer:DeligneCohomology})
and it is originally due to Milgram.  One of the advantages of
this over the usual Milnor construction is that if $G$ is an
abelian group then so is the geometric bar construction, so the
process can be iterated.

\subsection{Equivalence of groupoids}  Many of the groupoids
under consideration, such as $\GG(M)$ in the introduction, depend
on certain choices and are really only defined up to equivalence,
and in fact should really be thought of as representing a stack.
Equivalence of groupoids translates to homotopy equivalence of the
corresponding classifying spaces.  This means that the
constructions here ought to be natural with respect to such
equivalences.  This does have the advantage that certain
properties can just be proved for groups as it is shown below that
groupoids are equivalent to the disjoint union of groups.

\subsubsection{Natural transformations and homotopies}
 If $F$ and $F'$ are both functors from
$\cG$ to $\cH$, then a natural transformation, $T$, from $F$ to
$F'$, written $T\colon F\nattrans F'$ consists of a morphism
$T(x)\colon F(x)\to F'(x)$ for every object $x$ such that for
every morphism $g\colon x_1\to x_2$ in $\cG$ we have $F'(g)\circ
T(x_1)=T(x_2)\circ F(g)$, so $T$ \emph{intertwines} the two
functors. Two functors between groupoids are said to be
\emph{isomorphic} if there is a natural transformation between
them.

\begin{figure}[ht]
$$
{\relabelbox \small \epsfbox{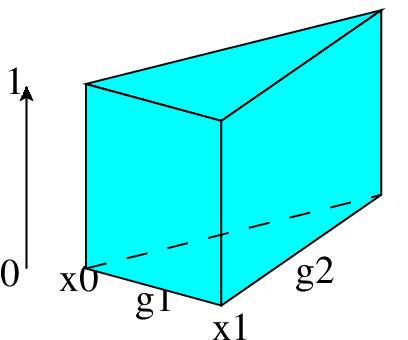} \special{ps: 0.14 0.14
scale} \relabel {g1}{$g_1$}
 \relabel {g2}{$g_2$}
 \relabel {x1}{$x_1$}
 \relabel {x0}{$x_0$}
 \relabel {x2}{$x_2$}
\relabel{1}{$1$} \relabel{0}{$0$}
 \endrelabelbox}
\longrightarrow {\relabelbox \small\epsfbox{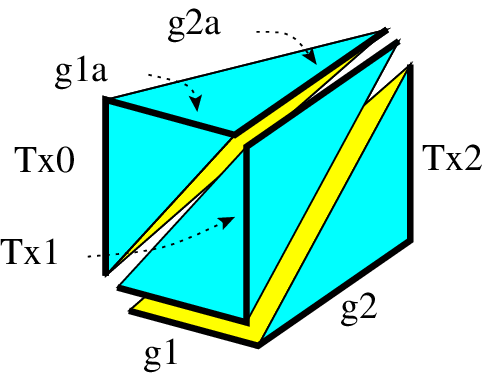}
\relabel{g1}{$F(g_1)$} \relabel{g1a}{$F'(g_1)$}
\relabel{g2}{$F(g_2)$} \relabel{g2a}{$F'(g_2)$}
\relabel{Tx0}{$T(x_0)$} \relabel{Tx1}{$T(x_1)$}
\relabel{Tx2}{$T(x_2)$}
\endrelabelbox}$$
\caption{The homotopy construction.} \label{Figure:Homotopy}
\end{figure}

Given such a natural transformation $T\colon F\nattrans F'$ we get
a homotopy $H_T$ between the classifying maps $\cB F$ and $\cB
F'$, in other words we get a map $H_T\colon \cB\cG\times [0,1]\to
\cB\cH$ such that $H_T({-},0)=F({-})$ and $H_T({-},1)=F'(-)$.
Explicitly it can be written as
\begin{multline*}
  H_T([g_n,t_n\;|\ldots|\;g_1,t_1],t)\\
    :=[F'(g_n),t_n\;|\ldots|\;F'(g_{i+1}),t_{i+1}\;|\;T({x_i}),t\;
          |~F(g_i),t_i\;|\ldots|\;F(g_1),t_1]\\
    \text{for }0\le t_n\le\dots\le t_{i+1}\le t\le t_i\le
    \dots t_1\le 1.
\end{multline*}
This formula should be made clearer by looking at
Figure~\ref{Figure:Homotopy}, which is where the formula comes
from.

\subsubsection{Equivalence of groupoids to unions of groups}
\label{Section:Equivalence}
 Two groupoids $\cG$ and $\cH$ are said
to be \emph{equivalent} if there are functors $F\colon \cG\to \cH$
and $\overline F\colon \cH\to\cG$ such that $F\circ\overline F$ is
isomorphic to the identity functor on $\cH$ and $\overline F\circ
F$ is isomorphic to the identity functor on $\cH$.  Equivalent
groupoids share many features, such as having the same
representation theory, having homotopic classifying spaces ($\cB
F$ and $\cB\overline F$ supplying the homotopy equivalence), and
having isomorphic homology and cohomology groups.

A fundamental feature of groupoids is that every groupoid is
equivalent to a disjoint union of groups, thus as we will see,
many problems about groupoids reduce to problems about groups. As
this will be so useful, I will give the details here.

To see this equivalence, it suffices to consider connected
groupoids and to show that every connected groupoid is equivalent
to a groupoid with one object, ie.\ a group.  So given a
\emph{connected} groupoid $\cG$ fix an object $x$, and define the
automorphism group $\Aut(x)$ of $x$ to be the group of morphisms
$\Mor_\cG(x,x)$.  Let $\overline {\Aut(x)}$ be the groupoid with
the single object, $x$, and the morphism group $\Aut(x)$.  Define
$i\colon \overline {\Aut(x)} \to \cG$ to be the inclusion.  To
construct an inverse equivalence to $i$ we need to make some
choices: for each object $y$ in $\cG$ pick a morphism $f_y\colon
y\to x$, taking $f_x\colon x\to x$ to be the identity morphism.
Now define the ``retract'' $r\colon \cG\to \overline {\Aut(x)}$ by
\begin{align*}
r(y)&:=x \qquad \text{for every object }y\\
r(g)&:=f_z\circ g\circ f_y^{-1} \quad\text{for every morphism }
        g\colon y\to z.
\end{align*}
To see that this gives an equivalence of groupoids, we just note
that $r\circ i=\id$ and that a natural isomorphism $T\colon
\id_\cG\nattrans i\circ r$ is given by defining $T(y):=f_y$ for
every object $y$.

From this the disconnected case easily follows:
\begin{thm}
If $\cG$ is a groupoid then $\cG$ is equivalent to the groupoid
$\coprod_x \overline{\Aut(x)}$ where $x$ runs over one object from
each connected component of $\cG$.
\end{thm}


\subsection{Loop groupoids, the Parmesan Theorem and transgression}
\label{Section:LoopGroupoidsParmesan}
 This is the heart of this
section on classifying spaces. Firstly in
Section~\ref{Subsubsection:LoopGroupoids} given a finite groupoid
$\cG$ the key notion of its loop groupoid $\LG$ is introduced.
Then in Section~\ref{Section:Parmesan} the Parmesan construction
is described, this is a homotopy equivalence $\cB\LG\to \cL\cB\cG$
where $\cL\cB\cG$ is the free loop space on the classifying space
$\cB\cG$.  In Section~\ref{Section:Transgression} this leads to
the $S^1$--transgression map on the level of cocycles $\tau\colon
Z^*(\cG,U(1))\to Z^{*-1}(\Lambda\cG,U(1))$, which will play a
leading role in Section~\ref{Section:TwistedStuff}.

\subsubsection{Loop groupoids}\label{Subsubsection:LoopGroupoids}
\begin{figure}[ht]
$$\vcenter{\epsfbox{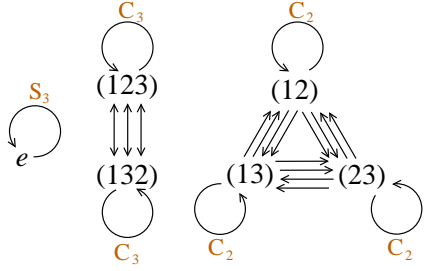}}$$
\caption{The loop groupoid $\Lambda\overline {S_3}$ or
equivalently the action groupoid $\cG_{S_3}(S_3^c)$ for $S_3$ with
the conjugation action (where $S_3$ denotes the symmetric group on
three letters). The automorphism groups have been marked on.}
\label{Figure:ActionCat}
\end{figure}

 For a finite groupoid $\cG$ define the loop groupoid (also
known as the inertia groupoid) $\Lambda\cG$ as follows. The
objects of $\LG$ are the self morphisms of $\cG$, ie.\
$\amalg_{x\in Ob\cG} \Mor_\cG(x,x)$, and there is a morphism
$g_\gamma\colon \gamma\to g\gamma g^{-1}$ in $\Lambda\cG$ for each
$\gamma\in\Mor_\cG(x,x)$ and $g\in \Mor_\cG(x,y)$ --- by abuse of
notation I will usually denote $g_\gamma$ it by $g$ when the
source is clear. Composition is induced from that in $\cG$.  An
example is pictured in Figure~\ref{Figure:ActionCat}.

An equivalent way to define $\LG$, and the point of view taken in
the introduction, is as $\Fun(\overline \Z,\cG)$, the functor
category to $\cG$ from the one object category with morphism group
$\Z$, with the morphisms in the functor category being given by
natural transformations.

As the name suggests this is supposed to be a model for ``loops on
$\cG$'', actually it is a model for loops on $\cB\cG$, as shown
below.  This should appear reasonable from the second description
in terms of the functor category as it is asserting that there is
a homotopy equivalence $\cB\Fun(\overline \Z,\cG)\simeq
\Maps(S^1,\cB\cG)$ and it is a standard fact that $B\overline \Z$
is homotopy equivalent to a circle.  This homotopy equivalence was
established by Neil Strickland \cite{Strickland}, but we will need
an explicit form of the map here.

\subsubsection{The Parmesan Theorem}
\label{Section:Parmesan} 
\begin{figure}[hbt]
$$
 \vcenter{\relabelbox\small\epsfbox{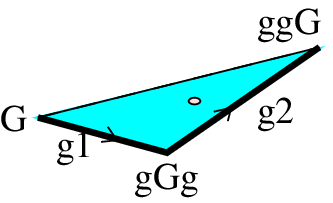}
 \relabel{G}{$\gamma$}
 \relabel{g1}{$g_1$}
 \relabel{g2}{$g_2$}
 \relabel{gGg}{${}^{g_1}\gamma$}
 \relabel{ggG}{${}^{g_2g_1}\gamma$}
 \endrelabelbox}
 \rightsquigarrow
 \vcenter{\relabelbox\small\epsfbox{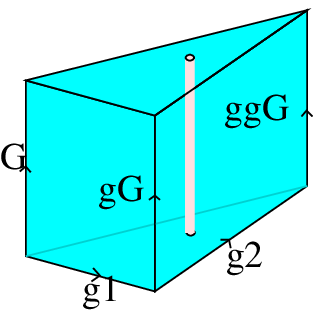}
 \relabel{G}{$\gamma$}
 \relabel{g1}{$g_1$}
 \relabel{g2}{$g_2$}
 \relabel{gG}{${}^{g_1}\gamma$}
 \relabel{ggG}{${}^{g_2g_1}\gamma$}
 \endrelabelbox}
\rightsquigarrow
 \vcenter{\relabelbox\small \epsfbox{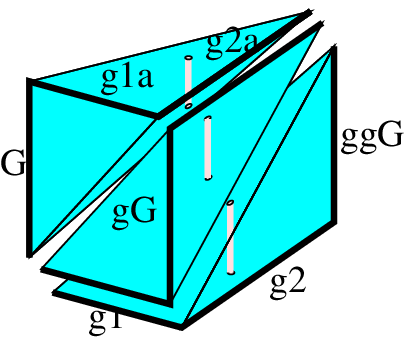}
 \relabel{G}{$\gamma$} \relabel{g1}{$g_1$} \relabel{g1a}{$g_1$}
 \relabel{g2}{$g_2$} \relabel{g2a}{$g_2$}
 \relabel{ggG}{${}^{g_2g_1}\gamma$}
 \relabel{gG}{${}^{g_1}\gamma$}
\endrelabelbox}
$$
\caption{The parmesan construction mapping a point in a
two-simplex of $\cB\LG$ to a loop in $\cB\cG$.  For reasons of
space, here ${}^g\gamma$ denotes $g\gamma g^{-1}$.}
\label{Figure:Parmesan}
\end{figure}

For simplicity of notation I will denote
the $n$--simplex in $\LG$
 $$g_n\ldots g_1\gamma g_1^{-1}\ldots g_n^{-1}
 \stackrel{g_n}\leftarrow\dots\stackrel{g_1}\leftarrow \gamma
 \qquad\text{ by }\qquad
 [g_n|\dots|g_1]\gamma.$$
We now get to the key construction. Denote by $\cL\cB\cG$ the free
loop space on $\cB\cG$, in other words
$\cL\cB\cG:=\Maps(S^1,\cB\cG)$, and define the Parmesan map
$\Par\colon\cB\Lambda\cG\to \cL\cB \cG$ by
\begin{multline*}
  \Par([g_n,t_n| \ldots|g_1,t_1]\gamma)(t)\\
  =[\;g_n,t_n\;|\ldots|\;g_{i+1},t_{i+1}\;|\;
              (g_i\ldots g_1)\gamma(g_i\ldots g_1)^{-1},t\;|\;
                     g_{i},t_{i}\;|\ldots|\;g_1,t_1\;]\\
  \text{ for\ }t_i\ge t \ge t_{i+1}.
\end{multline*}
This is illustrated graphically in Figure~\ref{Figure:Parmesan},
which is where it gets its name from.

Note that it was precisely to get this simple formula that I used
these coordinates on the simplices.  The key result about this map
is the following.
\begin{thm}\label{Thm:GroupoidParmesan}
If $\cG$ is a finite groupoid and $\Lambda \cG$ is its loop
groupoid then the Parmesan map $\Par\colon\cB\Lambda\cG\to \cL\cB
\cG$ is a homotopy equivalence.
\end{thm}
\begin{proof}
The main point in the proof is that we first prove it for finite
groups and then we prove it for general finite groupoids.  The
statement of the theorem for finite groups is
Theorem~\ref{Thm:GroupParemesan}, which is proved below.  So we
can assume it is true for a finite group.  Now assume that $\cG$
is a connected, finite groupoid --- the result in the
non-connected case follows from this. Pick any object $x$ in $\cG$
and let $G$ be the automorphism group of $x$. By
Section~\ref{Section:Equivalence}, the inclusion $\oG\to G$ gives
equivalence of groupoids, so the induced inclusion
$\cB\oG\to\cB\cG$ is a homotopy equivalence, and as the loop space
functor is a homotopy functor, the inclusion $\cL\cB\oG\to
\cL\cB\cG$ is a homotopy equivalence.  Similarly the inclusion
$\Lambda \oG\to \Lambda\cG$ is an equivalence, implying that the
inclusion $\cB\Lambda\oG\to \cB\Lambda\cG$ is a homotopy
equivalence.

It is immediate that the following diagram commutes
$$\begin{array}{ccc}
\cB\Lambda\cG&\stackrel\Par\longrightarrow &\cL\cB \cG\\
\llap{$\sim$}\uparrow&&\uparrow\rlap{$\sim$}\\
 \cB\Lambda\oG&\stackrel\sim\longrightarrow&\cL\cB\oG
\end{array},$$
and as the three marked maps are homotopy equivalences, it follows
that the Parmesan map must be as well.
\end{proof}

\subsubsection{The transgression map}
\label{Section:Transgression}
 As we have the map $\cB\Lambda\cG\to
\cL\cB\cG$ we also have the adjoint map (which can be thought of
as the evaluation map) $\cB\Lambda\cG\times S^1\to \cB\cG$ and
thus the map on the level of chains $C_n(\Lambda\cG\Z)\to
C_{n+1}(\cG,\Z)$ which is given by
\[
  [g_n| \ldots|g_1]\gamma\; \mapsto
  \sum_{i=0}^n (-1)^{n-i-1}[\;g_n\;|\ldots|\;g_{i+1}\;|\;
              (g_i\ldots g_1)\gamma(g_i\ldots g_1)^{-1}\;|\;
                     g_{i}\;|\ldots|\;g_1\;].
\]
Apart from the signs, this is seen for the case $n=2$ in
Figure~\ref{Figure:Parmesan}.

Dualizing this and taking $U(1)$--coefficients, we get the
transgression map on the level of cocycles:
\begin{thm}
For a finite groupoid $\cG$ the $S^1$--transgression map on the
cocycle level, $\tau\colon Z^{*}(\cG,U(1))\to Z^{*-1}(\LG,U(1))$,
is given on a cocycle $\alpha\in Z^{n}(\cG,U(1))$ by
\begin{multline*}
 \tau(\alpha)([g_{n-1}| \ldots|g_1]\gamma):=\\
    \prod_{i=0}^{n-1}
        \alpha([\;g_{n-1}\;|\ldots|\;g_{i+1}\;|\;
              (g_i\ldots g_1)\gamma(g_i\ldots g_1)^{-1}\;|\;
                     g_{i}\;|\ldots|\;g_1\;])^{(-1)^{n-1-i}}.
\end{multline*}
\end{thm}
We will use this for one-cocycles, two-cocycles and
three-cocycles, so it is worth explicitly recording them here.  So
supposing that $\omega$ is a three-cocycle, $\theta$ is a
two-cocycle and $\alpha$ is a one-cocycle, then
\[
  \tau(\alpha)([]\gamma):=\alpha([\gamma]);\qquad
  \tau(\theta)([g]\gamma)
     :=\frac{\theta([g\gamma g^{-1}|g])}{\theta([g|\gamma])};
\]
\[
  \tau(\omega)([h|g]\gamma)
     :=\frac{\omega([h|g|\gamma])
          \omega([hg\gamma g^{-1}h^{-1}|h|g])}
          {\omega([h|g\gamma g^{-1}|g])}
\]
The latter formula should look familiar to the reader \textit{au
fait} with the twisted Drinfeld double of a finite group.

\subsection{Action groupoids and the Parmesan Theorem for finite groups}
\label{Section:ParmesanForGroups}
\subsubsection{The definition of action groupoids}
\label{Subsection:ActionGroupoid}
Suppose that $G$ is a finite
group and $X$ is a finite $G$--set, so that there is an action
$G\times X\to X$, $(g,x)\mapsto g\cdot x$.  The \emph{action
groupoid} $\cG_G(X)$, which is sometimes denoted $G\times
X\rightrightarrows X$, is the groupoid whose objects are elements
of $X$, whose morphisms are of the form $g\cdot x \stackrel
g\leftarrow x$, with $x\in X$ and $g\in G$, and where composition
is given by the obvious
\[
  ( hg\cdot x \stackrel{h}{\leftarrow}g\cdot x)
         \circ ( g\cdot x \stackrel{g}{\leftarrow} x)
      :=( hg\cdot x \stackrel{hg}{\leftarrow} x).
\]
We will see below that the invariants of this groupoid are well
known things: its classifying space is the Borel construction on
$X$, its homology is the equivariant homology of $X$ and its
representation group is the equivariant K-theory of $X$.  Of
course as $X$ is a finite set, these things are not as exciting as
they might be, but we do get interesting information as we will
discover.

\subsubsection{Examples of action groupoids}
For any finite group $G$ there are three
finite $G$--sets that come to mind.

Firstly there is the trivial, one-element $G$--set $\{*\}$.  In
this case the action groupoid $\cG_G(\{*\})$ is just $\oG$.  Its
classifying space is $\cB G$, the classifying space of the group.

Secondly there is the $G$--set $G^L$, this is $G$ with the left
regular action on itself.  In this case $\cG_G(G^L)$ has $G$ as
its set of objects and from $g$ to $h$ there is a unique morphism
$g\stackrel{hg^{-1}}\longrightarrow h$.  In particular this means
that $\cG_G(G^L)$ is connected and has trivial automorphism
groups, so is equivalent to the trivial groupoid with one object,
thus the classifying space is contractible.  Furthermore $G$ acts
freely on the right on the groupoid $\cG_G(G^L)$ via
$(g\stackrel{hg^{-1}}\longrightarrow h)\cdot k:=
gk\stackrel{hg^{-1}}\longrightarrow hk$.  From this we get a free
action of $G$ on the classifying space.  We will denote the
classifying space $\cB\cG_G(G^L)$ by $\cE G$.  The quotient map by
the $G$--action on this classifying space gives a $G$--bundle $\cE
G\to\cB G$ known as the universal $G$--bundle, it is the
classifying map of groupoids coming from the $G$--set map $G^L\to
\{*\}$.

Thirdly there is the $G$--set $\Gconj$, this is $G$ with the
\emph{conjugation} action on itself, ie.\ $G\times \Gconj \to
\Gconj$ is $(g,h)\to ghg^{-1}$.  Unpacking the definition of the
action groupoid $\cG_G(\Gconj)$ reveals that it is precisely
$\Lambda \oG$, the loop groupoid on the one object category $\oG$.
The main result of this section is to show that the classifying
space of this is, via the Parmesan map, homotopy equivalent to the
loop space on the classifying space $\cB G$.

\subsubsection{Classifying spaces and the transgression map}
For a general $G$--set $X$ the classifying space of the action
groupoid $\cG_G(X)$ can be identified with $\cE G\times_G X$, the
Borel construction on $X$, as follows.  Firstly, to simplify
notation, I will denote a general $n$-simplex of $\cG_G(X)$ of the
form
 \[g_n\ldots g_1 \cdot x\stackrel {g_n}\leftarrow \dots
   \stackrel {g_2}\leftarrow g_1\cdot x \stackrel {g_1}\leftarrow x
   \quad\text{ by }\quad
   [g_n|\dots |g_1]x.\]
Recalling that the $\cE G$ was defined above to be the classifying
space $\cB \cG_G(G^L)$, we can define a map
\[
   \cE G\times_G X\to \cB \cG_G(X);\quad
  ([g_n\;|\dots |\;g_1]g,x)\mapsto [g_n\;|\dots |\;g_1]\;(g\cdot x).
\]
This is clearly a homeomorphism as an inverse is given by
\[
 [g_n|\dots |g_1]x\mapsto ([g_n|\dots |g_1]e,x).
\]
where $e$ is the identity in $G$.  Thus the classifying space on
the action groupoid of $X$ is the Borel construction on $X$.

Two facts worth mentioning here are: firstly, that the equivariant
cohomology of $X$ is defined to be the cohomology of the Borel
construction, this means that there is an isomorphism
$H^*_G(X)\cong H^*(\cG_G(X))$; and secondly, that the equivariant
K-theory of $X$ is by definition the Grothendieck group of
equivariant bundles on $X$, but an equivariant bundle on $X$
translates into a representation of the action groupoid of $X$,
thus $K_G(X)\cong \Rep(\cG_G(X))$, and this is not, in general the
same thing as the K-theory of the Borel construction.

 We are now ready to prove the
Parmesan Theorem for a finite group.  (Recall that $\cL\cB G$
denotes the free loop space on the classifying space of $G$.)

\begin{thm}
\label{Thm:GroupParemesan}%
If $G$ is a finite group then the Parmesan map
$\Par\colon\cB\Lambda\oG\to\cL\cB G$ (constructed in
Section~\ref{Section:Parmesan}) is a homotopy equivalence.
\end{thm}
\begin{proof}
We have already observed that the loop groupoid $\Lambda \oG$ for
a finite group is the same as the action groupoid $\cG_G(\Gconj)$
for $G$ with the conjugation action, and we have an explicit
homeomorphism $\cB\cG_G(\Gconj)\cong \cE G\times_G \Gconj$ so we
will consider the Parmesan map as a map $\cE G\times_G \Gconj\to
\cL\cB G$.  We will show that this is a homotopy equivalence by
constructing a homotopy inverse.

Define $\varphi\colon\cL\cB\oG\to \cE G\times_G \Gconj$ as
follows. If $\gamma$ is a loop in $\cL\cB\oG$, thought of as a map
$[0,1]\to \cB\oG$ with the endpoints mapping to the same point,
then pick a lift of $\gamma(0)$ to the covering space $\cE G$.
This leads to a lift $\tilde\gamma\colon [0,1]\to \cE G$ of
$\gamma$, and as $\tilde\gamma(0)$ and $\tilde\gamma(1)$ lie in
the same fibre, they differ by an element $g_{\tilde\gamma}\in G$,
so $(\tilde\gamma(0),g_{\tilde\gamma})$ gives an element in $\cE
G\times G$, but this element depended on the lift of $\gamma(0)$,
however this ambiguity is precisely absorbed by considering this
as an element $\varphi(\gamma)$ of $\cE G\times_G\Gconj$. A proof
that this map is a homotopy equivalence is given in \cite[Lemma
2.12.1]{Benson:BookII}.

A two-line calculation shows that $\varphi\circ\Par$ is the
identity.  As $\varphi$ is a homotopy equivalence, it follows that
$\Par$ is too.
\end{proof}

Thus in the case of finite group $G$ we get an explicit map
$\cB\cG(\Gconj)\to\cL\cB G$.  We know that a groupoid is
equivalent to the disjoint union of its automorphism groups and in
this case the automorphism group of $x\in \Gconj$ is just its
centralizer $C_x$ in $G$, that is, the set of elements in $G$
which commute with $x$.  Thus we get the standard result that
$H^*(\cL\cB G,U(1))\cong \bigoplus H_\text{gp}^*(C_x,U(1))$ where
$x$ runs over one representative for each conjugacy class.
Combining this with the transgression map
$H_{\text{gp}}^*(G,U(1))\to H^{*-1}(\cL\cB G,U(1))$ gives maps
$\tau_x\colon H^*_\text{gp}(G,U(1))\to
H^{*-1}_\text{gp}(C_x,U(1))$ for any element $x\in G$.  As these
maps crop up in a few places such as discrete torsion, it is worth
writing them down here on the level of cocycles in small degree.
If if $\theta\in Z^2_\text{gp}(G,U(1))$ is a group two-cocycle on
$G$, $\omega\in Z^3_\text{gp}(G,U(1))$ is a three-cocycle, $x\in
G$ and $h,g\in C_x$ then

\[
  \tau_x\theta(g):=
        \frac{\theta(g,x)}{\theta(x,g)};\qquad
  \tau_x\omega(h,g):=
        \frac{\omega(h,g,x)\omega(x,h,g)}{\omega(h,x,g)}.
\]

\section{Twisted representations and cohomology}\label{Section:TwistedStuff}
This section is the algebraic heart of the paper in which the
theory of twisted representations of finite groupoids is
developed.  The philosophy behind it is very much that of
$n$--gerbes.  What this means is that starting with $n=0$ and
working up to $n=2$ we see how an $n$--cocycle leads to an
$(n-1)$--category which is to be interpreted as the integral or
space of sections or push-forward to a point of the `complexified
gadget' associated to the $(n-1)$--gerbe. Furthermore the things
at level $(n-1)$ are classified in an appropriate sense by the
transgressed $(n-2)$ cocycle over the loop groupoid.  All the
constructions are natural with respect to equivalences of
groupoids.

First here is a quick reminder on the theory of twisted
representations of finite groups --- a reasonable reference is
\cite{Karpilovsky:Book}.

\subsection*{Reminder on twisted representations of finite groups}
Suppose that $G$ is a finite group and $\theta\in Z^2(G,U(1))$ is
a normalized two-cocycle, so it is a $U(1)$--valued function on
$G\times G$ satisfying:
\begin{align*}
\theta(h,g)\theta(kh,g)^{-1}\theta(k,hg)\theta(k,h)^{-1}&=1
     \qquad\text{for all }k,h,g\in G;\\
\theta(e,g)=\theta(g,e)&=1
     \qquad\text{for all }g\in G.
\end{align*}
A \emph{$\theta$--twisted representation} of $G$ (also known as a
projective representation with twisting-cocycle $\theta$) is
defined to be a function $\rho\colon G\to \End(V)$ for some vector
space $V$ such that
\begin{align*}
\rho(h)\rho(g)&=\theta(h,g)\rho(hg);\\
\rho(e)&=1.
\end{align*}
There are a couple of other ways to think of a twisted
representation.  One way is via the central extension $U(1)\to
\pretheta G\to G$, where the underlying set of $\pretheta G$ is
taken as $U(1)\times G$ and the product is taken to be
$(z,h)(w,g):=(\theta(h,g)zw,hg)$; a $\theta$--twisted
representation of $G$ can then be thought of as a representation
of $\pretheta G$ in which the central $U(1)$ acts in its natural
way.  Another way to think about a $\theta$--twisted
representation is as a representation of the twisted group algebra
${}^\theta\!\C G$, which is the algebra with basis elements
$\{\langle g\rangle\}_{g\in G}$ and with multiplication given by
$\langle h\rangle \langle g\rangle:=\theta(h,g)\langle hg\rangle$.
To reconcile these two other points of view, note that
${}^\theta\!\C G\cong\C\otimes_{U(1)} {}^\theta\!G$.

One of the simplest non-trivial examples is the following.
Consider the four-group $\Z_2\times \Z_2=\langle a,b\mid
a^2=b^2=(ab)^2=1\rangle$, and define the two-cocycle $\theta_V\in
Z^2(\Z_2\times \Z_2,U(1))$ by
\[\theta(a^{u_a}b^{u_b},a^{v_a}b^{v_b})=(-1)^{u_a.v_b}.\]
Then the map $\pretheta\C(\Z_2\times \Z_2)\to M_2(\C)$ given by
\[
  e\mapsto \left( \begin{array}{cc}1&0\\0&1\end{array}\right),
  \quad
  a\mapsto \left( \begin{array}{cc}1&0\\0&-1\end{array}\right),
  \quad
  b\mapsto \left( \begin{array}{cc}0&i\\-i&0\end{array}\right),
  \quad
  ab\mapsto \left( \begin{array}{cc}0&-i\\-i&0\end{array}\right),
\]
is an isomorphism.

 It is well-known that ordinary representations
of finite groups are classified up to equivalence by their
character, this is the function on $G$, invariant under
conjugation, given on an element by taking the trace of the
element in the corresponding representation.  Trying the exact
same trick with twisted representations does not work in general
as taking the trace does not lead to a function on $G$ invariant
under conjugation, as the following lemma demonstrates.

\begin{lemma}
\label{Lemma:TwistedConjugation}
 For $\rho$ a $\theta$--twisted
representation of the finite group $G$, with $g,h\in G$, the
following conjugation relation holds:
\[
   \Tr_\rho(hgh^{-1})
     =\frac{\theta(hgh^{-1},h)}{\theta(h,g)}\Tr_\rho(g).
\]
\end{lemma}
\begin{proof}
Firstly using the basic property of $\rho$ twice reveals
\[
  \rho(hgh^{-1})
    =\theta(hg,h^{-1})^{-1}\theta(h,g)^{-1}\rho(h)\rho(g)\rho(h^{-1}).
\]
Then applying it again to $\rho(hh^{-1})$ gives
$\rho(h^{-1})=\theta(h^{-1},h)\rho(h)^{-1}$.  Substituting this
into the above and using the two-cocycle relation
$\theta(hg,h^{-1})\theta(h^{-1},h)=\theta(hgh^{-1},g)$ coming from
the triple $(hg,h^{-1},h)$ gives
\[ \rho(hgh^{-1})
     =\theta(hgh^{-1},h)\theta(h,g)^{-1}\rho(h)\rho(g)\rho(h)^{-1}.
\]
Applying trace gives the requisite relation.
\end{proof}
We can define a $\theta$--twisted character to be a function on
$G$ which satisfies the conjugation relation in
Lemma~\ref{Lemma:TwistedConjugation}, and then you can show that
$\theta$--twisted characters classify $\theta$--twisted
representations up to equivalence (see \cite{Karpilovsky:Book}).
We will see below that $\theta$--twisted characters can be thought
of as sections of a certain line bundle over the loop groupoid
$\Lambda G$.

\subsection{Zero-cocycles, locally-constant functions and
integration} To a zero-cocycle on a groupoid we can associate its
integral which is just a complex number.

Suppose that $\cG$ is a finite groupoid. A zero-cycle $\beta\in
Z^0(\cG,U(1))$ is the same thing as a locally constant function,
that is a $U(1)$--valued function which is constant on connected
components of $\cG$.  We can define $\int_\cG \beta$ the integral
of $\beta$ as follows
\[
  \int_\cG\beta:=\sum_{x\in\Ob(\cG)} \frac{\beta(x)}{|x\to|},
\]
where $|x\to{}|$ denotes the number of morphisms in $\cG$ with $x$
as their source.  Essentially this amounts to picking a measure on
$\cG$ and the point about this measure is that it is invariant
under pull-back via equivalences, ie. if $F\colon \cH \to \cG$ is
an equivalence then $\int_\cG \beta=\int_\cH F^*\beta$.  This
follows from an alternative expression for the integral:
\[
  \int_\cG\beta=\sum_{[x]\subset\Ob(\cG)} \frac{\beta(x)}{|\Aut(x)|},
\]
where $x$ runs over one object from each connected component of
$\cG$.  This is immediate as $\beta$ is constant on connected
components and $|[x]||\Aut(x)|=|x\to{}|$.

\subsection{One-cocycles, circle-bundles and flat sections.}
To a one-cocycle on a groupoid we can associate a `flat complex
line bundle' and from that its `space of flat sections' which is a
complex vector space.
\subsubsection{Circle bundles and associated line bundles} For a
finite groupoid $\cG$ a one-cocycle $\alpha\in Z^1(\cG,U(1))$ is
the same thing as a functor $\cG\to \overline{U(1)}$, where as
usual $\overline{U(1)}$ means the one object groupoid with $U(1)$
as its group of morphisms --- the cocycle condition translates
into functoriality under composition of morphisms.  Such a functor
can also be thought of as a ``trivialized $U(1)$--bundle with flat
connection'' over $\cG$, with a copy of the standard $U(1)$
sitting over each object of $\cG$ and with parallel transport
along a morphism given by the cocycle evaluated on that morphism.

As an aside, note that more generally, a ``$U(1)$--bundle with
flat connection'' over $\cG$ would be a functor into the category
of principal $U(1)$--spaces (also known as $U(1)$--torsors).  This
category is equivalent to the one object category
$\overline{U(1)}$, reflecting the fact that in this finite
setting, all $U(1)$--bundles are trivializable.

Associated to $\alpha$ as above we have the corresponding
``trivialized complex-line bundle with flat connection'' which we
can denote $\alpha_\C\colon \cG\to\overline\C$, this has a
standard copy of $\C$ over each object and the same parallel
transport as the $U(1)$--bundle.  A \emph{flat section} of the
trivialized complex line bundle $\alpha_\C$ is then defined to be
a map $s\colon\Ob(\cG)\to \C$ such that for any morphism $f\colon
x_1\to x_2$ in $\cG$ we have $\alpha(f)s(x_1)=s(x_2)$.  It is then
a simple observation that for each connected component of $\cG$
there is either a one-dimensional or zero-dimensional space of
flat sections, depending whether or not $\alpha$ sends all
automorphism in that component to the identity in $U(1)$.  Below
we will see a formula for the number of flat sections of
$\alpha_\C$ in general.

\subsubsection{Inner product on the space of sections}
Suppose that $s_1$ is a section of the line bundle $\alpha_\C$
over the groupoid $\cG$, then this is a function
$s_1\colon\Ob(\cG)\to \C$ with certain transformation properties.
Now the complex conjugate $\overline s_1\colon\Ob(\cG)\to \C$ is
actually a section of $(\alpha^{-1})_\C$ so, given another flat
section $s_2$, the product $\overline {s_1}s_2\colon \Ob(\cG)\to
\C$ is just a flat section of the trivial line bundle, or in other
words, a function constant on connected components. Thus defining
an inner product
\[  \langle s_1,s_2\rangle:=\int_\cG\overline {s_1 }s_2\]
gives something that is invariant under equivalence of groupoids.

\subsubsection{Number of flat sections}
We can now see the formula for the number of flat sections of the
flat line bundle $\alpha_\C$ corresponding to the one-cocycle
$\alpha$.  This can be calculated using the transgression of
$\alpha$ which is a zero-cocycle on the loop groupoid (see
Section~\ref{Section:Transgression}).
\begin{thm} \label{Theorem:DimSectionsBundle}
If $\alpha\in Z^1(\cG,U(1))$ then the number of flat sections of
the flat line bundle $\alpha_\C$ is given by integrating
$\tau(\alpha)\in Z^0(\Lambda\cG,U(1))$, the transgression of
$\alpha$, over the loop groupoid $\Lambda\cG$, in other words
\[
  \dim \Gamma_{\cG}(\alpha_\C)
         =\int_{\Lambda\cG}\tau(\alpha).
\]
\end{thm}
\begin{proof}
It suffices to prove this in the case that $\cG$ is connected as
the number of sections is clearly additive under disjoint union.
Furthermore, as everything is invariant under equivalence of
groupoids, it suffices to prove it in the case when $\cG$ has just
one object, in other words, when $\cG=\oG$ for some finite group
$G$.  In this case $\alpha\in Z^1(\oG,(U(1))$ is just a group
homomorphism $G\to U(1)$, and the line bundle $\alpha_\C$ has a
section precisely when $\alpha$ is trivial.

Looking at the right hand side of the above formula, by definition
we have
\[
  \int_{\Lambda\oG}\tau(\alpha)=
             \sum_{x\in G}\frac{\alpha(x)}{|G|}.
\]
However, for any $h$ in $G$ we have
\[
\sum_x\frac{\alpha(x)}{|G|}\alpha(h)=\sum_x
\frac{\alpha(xh)}{|G|}=\sum_x \frac{\alpha(x)}{|G|}.
\]
 So $\sum_x\frac{\alpha(x)}{|G|}(\alpha(h)-1)=0$, thus
 \[
   \int_{\Lambda\oG}\tau(\alpha)
   =\left\{
   \begin{array}{ll} 1&\alpha\equiv 1\\0 &\text{otherwise}
   \end{array}
   \right\}
   = \dim\{\text{sections of }\alpha_\C\}.
\]
\end{proof}
The above has a generalization, as $\alpha_\C\colon\cG\to
\overline \C$ can be thought of as a one dimensional
representation of $\cG$.  If we let $F\colon \cG\to\Vect$ be an
ordinary representation of the groupoid $\cG$ then we also have,
$\Gamma_\cG(F)$ the space of flat section of $F$, which can be
thought of as the \emph{invariant part} of the representation $F$.
The character of $F$ is then a locally constant function on the
loop groupoid $\Lambda\cG$ and the same proof as above gives the
following result.
\begin{prop}\label{Prop:DimInvariants}
  For $F\colon \cG\to \Vect$ a representation of a finite groupoid
  $\cG$ the number of flat sections is given by integrating the
  character $\chi_F$ over the loop groupoid:
  \[
    \dim \Gamma_{\cG}(F)
         =\int_{\Lambda\cG}\chi_F.
  \]
\end{prop}

\subsection{Two-cocycles, gerbes and twisted representations.}
This subsection follows the same pattern as the previous one, but
one categorical level higher.  Firstly, given $\theta$ a
two-cocycle on a finite groupoid $\cG$ a $\theta$--twisted
representation is interpreted as a section of a bundle of
categories associated to the corresponding `gerbe' over $\cG$.
These sections naturally form a category.  Such a twisted
representation can also be viewed as a representation of the
twisted groupoid algebra $\pretheta\C\cG$, which is seen to be
semisimple.  After an appropriate notion of inner product on the
category of twisted representations is defined, the idea of
character is introduced, and that the representations are
classified by their character.

\subsubsection{Basics} If $\cG$ is a finite groupoid
and $\theta\in Z^2(\cG,U(1))$ is a two-cocycle, then by analogy
with the finite group case we can define a $\theta$--twisted
representation of $\cG$ to be something like a representation of
$\cG$ --- ie.\ a functor $\cG\to\Vect$
--- but which fails the composition rule in a manner controlled by
$\theta$.  So we define, F, a $\theta$--twisted representation of
$\cG$ to be something which associates to each object $x$ a vector
space $F(x)$ and to each morphism $g\colon x_1\to x_2$ a linear
map $F(g)\colon F(x_1)\to F(x_2)$, such that
$F(g_2)F(g_1)=\theta([g_2|g_1])F(g_2\circ g_1)$ whenever $g_2$ and
$g_1$ are composable.  I will denote a twisted representation of
$\cG$ by $F\colon\pretheta\cG\to\Vect$, this is slightly ambiguous
given the notation that follows, but this should not cause any
confusion.

As in the finite group case, there are several other ways of
thinking of twisted representations.  Firstly, and philosophically
the gerbel way to think of them, is via the `trivialized central
extension' $\pretheta\cG$ defined as follows: objects of
$\pretheta\cG$ are the same as objects of $\cG$, but the morphism
sets are $\Mor_{\pretheta\cG}(x, y):=U(1)\times \Mor_\cG(x,y)$
with the twisted composition given by $(z_2,g_2)\circ
(z_1,g_1):=(\theta([g_2|g_1])z_2z_1,g_2\circ g_1)$. This is can be
thought of as a central extension written as
$\overline{U(1)}\to\pretheta\cG\to\cG$, which is what I will mean
by a gerbe (with flat connection) over $\cG$.  A $\theta$--twisted
representation of $\cG$ is then a representation of $\cG$ in which
the central $U(1)$ acts in its natural way.  An even more
appropriate way to think is via the associated `two-vector bundle'
$\Vect\to \pretheta\cG_{\Vect}\to \cG$ where we define
$\pretheta\cG_{\Vect}:=\Vect\times_{\overline{U(1)}}\pretheta\cG$,
so that%
\begin{align*}
   \Ob(\pretheta\cG_{\Vect})&:=\Ob(\Vect)\times \Ob(\cG),\\
   \Mor_{\pretheta\cG_{\Vect}}((V_1,x_1),(V_2,x_2))
      &:=\Mor_{\Vect}(V_1,V_2)\times\Mor_\cG(x_1,x_2),\\
   (\varphi_2,g_2)\circ(\varphi_1,g_1)
      &:=(\theta([g_2|g_1])\varphi_2\circ\varphi_1,g_2\circ g_1).
\end{align*}
A $\theta$--twisted representation of $\cG$ is then just a section
of $\pretheta\cG_{\Vect}\stackrel\pi\longrightarrow \cG$, in other
words a functor
$\cG\stackrel{F}\longrightarrow\pretheta\cG_{\Vect}$ such that
$\pi\circ F=\id_\cG$.

The collection of $\theta$--twisted representations of $\cG$ form
a category which I will denote by $\thetacRep(\cG)$.  The
morphisms are the intertwiners, so if $F,F'\colon \pretheta\cG \to
\Vect$ are representations then a morphism $T\colon F\nattrans F'$
is essentially a natural transformation, so for each $x\in
\Ob(\cG)$ there is a linear map $T(x)\colon F(x)\to F'(x)$ such
that for any $g\colon x\to y$ we have $T(y)F(g)=F(g)T(x)$.


Another way to think of a $\theta$--twisted representation of
$\cG$, which is relevant for the twisted Drinfeld double, is as an
ordinary representation of $\pretheta\C\cG$ the twisted groupoid
algebra, this is defined to be the algebra spanned by the
morphisms in $\cG$, such that the product $\langle
g_2\rangle\langle g_1\rangle$ is zero if $g_2$ and $g_1$ are not
composable in $\cG$ and is $\theta([g_2|g_1])\langle{g_2\circ
g_1}\rangle$ otherwise. Ordinary representations of this algebra
are then essentially the same thing as $\theta$--twisted
representations of the original groupoid, as is explained in the
next proposition.

\begin{prop} \label{Prop:AlgGpdRep}
If $\cG$ is a finite groupoid and $\theta\in Z^2(\cG,U(1))$ is a
two-cocycle then the category of representations of the twisted
groupoid algebra $\pretheta\C\cG$ is equivalent to the category
$\thetacRep(\cG)$ of $\theta$--twisted representations of the
groupoid $\cG$.
\end{prop}
\begin{proof}
Suppose that $F\colon \pretheta\cG\to \Vect$ is a
$\theta$--twisted representation of $\cG$, then define the vector
space $V:=\bigoplus_{x\in \Ob(\cG)}F(x)$ and the map $\rho\colon
\pretheta\C\cG\to \End(V)$ when $g\colon x_0\to x_1$ is a morphism
in $\cG$ by $\rho(\langle g\rangle)v:=F(g)v$ if $v\in F(x_0)$ and
zero if $v$ is in any other homogeneous component.

Conversely, suppose that $\rho\colon \pretheta\C\cG\to \End(V)$ is
a representation of the groupoid algebra.  To create a
$\theta$--twisted representation of the groupoid $\cG$ we have to
break the vector space $V$ into pieces: for $x\in Ob(\cG)$ define
the vector space $F(x):=\rho(\langle\id_x\rangle)V$.  The fact
that the identity in the twisted groupoid algebra splits into a
sum of the individual identities in the groupoid which are
indecomposable idempotents gives that $V$ is just the direct sum
$\bigoplus_{x\in \Ob(\cG)}F(x)$.  The functor
$F\colon\pretheta\cG\to \Vect$ is then completed by defining it on
morphisms as follows, if $g\colon x_0\to x_1$ is a morphism in
$\cG$ then $F(g):=\rho(\langle g\rangle)|_{F(x_0)}$.
\end{proof}

  As in the case of
ordinary representations of finite groups, every twisted
representation is completely reducible so can be written as a
direct sum of irreducible representations.
\begin{thm}
  If $F'$ is a subrepresentation of the $\theta$--twisted
  representation $F\colon \pretheta\cG\to \Vect$, then there is a
  complementary subrepresentation $F''$ so that $F=F'\oplus F''$.
\end{thm}
\begin{proof}
  This is similar to the usual proof for untwisted representations
  of groups, but takes a little more care.

  For each object $x$ in $\cG$ pick a projection $P(x)\colon
  F(x)\to F(x)$ so that $\Image(P(x))=F'(x)$, this means that
  $P(x)^2=P(x)$ and, as $F'$ is a subrepresentation, for any
  morphism $g\colon x\to y$ we have $P(y)F(g)P(x)=F(g)P(x)$.
  Typically, $P$ will not be an intertwiner, but in the usual way
  we can make it into one.  Define
   \[ P'(x):=\frac 1{|x{\to}|}\sum_{x\stackrel g \to y}
              F(g)^{-1}P(y)F(g),
   \]
  where the sum is over all morphisms coming out of $x$ and
  $|x{\to}|$ is the number of such morphisms.
  Then it is straight forward (but not without some grittiness)
  to show that $P'(x)^2=P'(x)$, that
  $\Image(P'(x))=F'(x)$ and that $P'$ is an intertwiner.   We can
  then take the complementary representation $F''$ to be $\ker
  P'$.
\end{proof}
\noindent It follows from this that the twisted groupoid algebra
$\pretheta\C\cG$ is semisimple.

 If we now denote the Grothendieck group of
$\thetacRep(\cG)$ the category of $\theta$--twisted
representations by $\thetaRep(\cG)$, then the above means that the
representation group $\thetaRep(\cG)$ is the free abelian group
generated by the equivalence classes of the irreducible
representations.

\subsubsection{Inner product on representations}
\label{Subsection:IPOnReps} We can define the functor
\[
  \langle {-},{-}\rangle\colon
   \thetacRep(\cG)\times\thetacRep(\cG)\to \Vect;
 \quad
 \langle F_1, F_2\rangle :=\Hom_{\thetacRep(\cG)}(F_1,F_2).
\]
 which is thought of as
an ``inner product'' on the category of $\theta$--twisted
representations. This sort of situation often arises in categories
in which there are internal hom functors and duals, this isn't the
case in the twisted situation and it is worth explaining this.

The first thing to do is to understand the contragredient twisted
representation.  If $F\colon \pretheta\cG\to \Vect$ is a
$\theta$--twisted representation then define $F^\vee$ to be the
``functor''
\begin{align*}
 F^\vee(x)&:=F(x)^\vee=\Hom_\C(F(x),\C);\\
  F^\vee(g\colon x_0\to x_1)
    &:=
    \bigl((F(g)^{-1})^\vee\colon F(x_0)^\vee\to F(x_1)^\vee\bigr).
\end{align*}
%
Note that in general $F(g)^{-1}$ is not the same as $F(g^{-1})$,
but they differ by a factor of $\theta(g,g^{-1})$. A two-line
calculation gives that $F^\vee$ is a $\theta^{-1}$--twisted
representation of $\cG$, whereas a one-line calculation gives that
the tensor product of a $\theta$--twisted representation and a
$\psi$--twisted representation is a $\theta\psi$--twisted
representation, so if $F_1$ and $F_2$ are both $\theta$--twisted
representations then $F_1^\vee\otimes F_2$ is a genuine
representation of $\cG$, in other words a genuine functor
$\cG\to\Vect$, so we can take the space of flat sections
$\Gamma_\cG(F_1^\vee\otimes F_2)$.  It is straight forward to see
that this is the hom-set $\Hom_{\thetacRep(\cG)}(F_1,F_2)$.

\subsubsection{Characters of twisted representations of groupoids}
\label{Section:TwistedCharacters}
We can now essentially copy the
construction from twisted representations of finite groups, and
classify, up to equivalence, $\theta$--twisted representations of
a finite groupoid $\cG$ by $\theta$--twisted characters, which are
flat sections of a line bundle over the loop groupoid
$\Lambda\cG$, the line bundle being $\tau(\theta)_\C$, the line
bundle associated to the circle-bundle transgressed from $\theta$.

Define the character map
  \[
    \chi\colon \thetaRep(\cG)\to
       \Gamma_{\Lambda\cG}(\tau(\theta)_\C);
    \quad
    F\mapsto \chi_F;
    \quad
    \chi_F(\gamma):=\Tr\circ F(\gamma).
  \]
The fact that $\chi_F$ is indeed a section of $\tau(\theta)_\C$,
ie.\ that
 \[
   \chi_F(h\gamma h^{-1})=
       \tau(\theta)([h]\gamma)\;
       \chi_F(\gamma),
 \]
is proved exactly as Lemma~\ref{Lemma:TwistedConjugation}.  We can
call the space of flat sections
$\Gamma_{\Lambda\cG}(\tau(\theta)_\C)$ the space of
$\theta$--twisted characters of $\cG$.  The fact that the twisted
character map $\chi$ is an injection follows from the following
result.
\begin{prop}
  The twisted character map respects the inner products: if
  $F_1$ and $F_2$ are $\theta$--twisted
  representations of $\cG$ then
  \[\langle \chi_{F_1},\chi_{F_2}\rangle =
     \dim\langle F_1,F_2\rangle.
  \]
\end{prop}
\begin{proof}
  There are two preliminary results needed here.  The first is
  that if $F$ and $F'$ are respectively $\theta$--twisted and $\phi$--twisted
  representations of $\cG$ then $\chi_{F\otimes F'}=
  \chi_F\cdot\chi_{F'}$ --- this follows immediately from the
  properties of traces.

  The second result needed is slightly more
  awkward to prove but it is that the character of a contragredient
  representation is the conjugate of the character of the original representation:
  $\chi_{F^\vee}=\overline{\chi_{F}}$.  Note firstly that the
  trace of the adjoint of a linear map is equal to the trace of
  the original map and if $A$ is a linear map with all eigenvalues
  in $U(1)$ then $\Tr A^{-1}=\overline{\Tr A}$.  So for $g$ any
  automorphism in $\cG$ we have
  $\chi_{F^\vee}(g)=\Tr((F(g)^{-1})^\vee)=\Tr(F(g)^{-1})$, thus if
  $F(g)$ has all its eigenvalues in $U(1)$ then this equals
  $\overline{\Tr(F(g))}=\overline{\chi_F(g)}$, and we are done.
  Hence it suffices to show that $F(g)$ has all its eigenvalues in
  $U(1)$.  Begin by supposing that the
  eigenvalues of $F(g)$ are $d_1,\dots,d_m$, and that $g$ has
  order $q$ in $G$.  The $F(g)^q$ differs from $F(g^q)=F(e)=\id$
  by a scalar which is a product of factors coming from $\theta$,
  and as $\theta$ is $U(1)$--valued, this scalar lives in $U(1)$.
  Thus the eigenvalues of $F(g)^q$ --- $d^q_1,\dots,d^q_m$ --- all
  live in $U(1)$ and hence the eigenvalues of $F(g)$ live in
  $U(1)$.

  We can now proceed with the proof of the proposition.  Using the above two
  results together with the ideas of Section~\ref{Subsection:IPOnReps} and
  Proposition~\ref{Prop:DimInvariants} we find
  \begin{align*}
    \dim\langle F_1,F_2\rangle &=
    \dim\Gamma_{\cG}(F_1^\vee\otimes F_2)=
    \int_{\Lambda\cG}\chi_{F_1^\vee\otimes F_2}=
    \int_{\Lambda\cG}\chi_{F_1^\vee}\cdot\chi_{F_2}\\
    &=\int_{\Lambda\cG}\overline{\chi_{F_1}}\cdot\chi_{F_2}
    =\langle \chi_{F_1},\chi_{F_2}\rangle.
  \end{align*}
%
\end{proof}
From this it is clear that $\theta$--twisted representations are
classified up to equivalence by their character, as the set of
(equivalence classes of) irreducible twisted representations is
mapped bijectively to an orthonormal set.  Even better is that
this is a full embedding so that the irreducible twisted
representations are actually mapped to an orthonormal basis.
\begin{thm}\label{Theorem:DimSectionsGerbe}
For $\cG$ a finite groupoid and $\theta\in Z^2(\cG,U(1))$, the
character map defined above induces an isomorphism of inner
product spaces:
\[\chi\colon \thetaRep(\cG)\otimes\C
\stackrel\simeq\longrightarrow
\Gamma_{\Lambda\cG}(\tau(\theta)_\C).
\]
\end{thm}
\begin{proof}
  This follows from Theorem~\ref{Theorem:DecomposeRep} below --- which
  decomposes $\thetaRep(\cG)$ as the representation groups of the
  automorphism groups of the components of $\cG$ --- together with
  the fact that the theorem holds for finite \emph{groups}
  \cite{Karpilovsky:Book}.
\end{proof}

Applying Theorem~\ref{Theorem:DimSectionsBundle} to get the
dimension of the above spaces gives a formula for the rank of the
twisted representation group as follows.
\begin{cor}\label{Cor:NumberTwistedReps}
For $\cG$ a finite groupoid, and $\theta\in Z^2(\cG,U(1))$,
  \[
     \#\{\text{irreducible $\theta$--twisted representations
     of $\cG$}\}=
          \int_{\Lambda^2\cG}\tau^2(\theta).
\]
\end{cor}
Moreover, specializing to the case when $\cG$ is just a group, we
can unpack the definition of the right-hand side to obtain:
\begin{cor} For $G$ a finite group and $\theta\in
Z^2_{\text{gp}}(G,U(1))$,
\[
      \#\{\text{irreducible $\theta$--twisted representations
     of $G$}\}=
          \frac{1}{|G|}\sum_{xg=gx}\frac{\theta(x,g)}{\theta(g,x)}.
\]
\end{cor}
\noindent This does not seem to appear in \cite{Karpilovsky:Book}
but was undoubtedly known to representation theorists: it is
equation (6.40) in \cite{DijkgraafWitten}.

The above formula can easily be applied to the example in the
introduction to this section, to verify that there is one
irreducible $\theta_V$--twisted representation of
$\Z_2\times\Z_2$, which is of course the two-dimensional
representation given.

\subsection{Classifying groupoid representations via group
representations.} From a groupoid point of view it is essentially
a triviality that a representation of a groupoid is determined by
its restriction to the automorphism groups, this is true in the
twisted case as well, but one has to be more careful there.

\subsubsection{Construction} \label{Subsubsection:DPRInduction}  We will
restrict for the moment to the case of connected finite groupoids.
Suppose that $x$ is an object of $\cG$ and that $\theta\in
Z^2(\cG,U(1))$ is a two-cocycle then there is a restriction
cocycle $\theta|_x\in Z^2(\Aut(x),U(1))$.  If
  $$\rho\colon{}^{\theta|_x}\!\Aut(x)\to \Aut(V)$$
is a $\theta|_x$--twisted representation then we can define a
$\theta$--twisted representation of $\cG$ in the following way.
First pick a map $f_y\colon y\to x$ for every object $y$, with
$f_x\colon x\to x$ being the identity.  Now define
\begin{alignat*}{2}
  &\text{for every object $y$ in $\cG$:}&\quad
  F_\rho(y)&:=V; \\
  &\text{for every $g\colon y\to z$:}&
  F_\rho(g)&:=\frac{\theta([f_z|g])}
              {\theta([f_z g f_y^{-1}|f_y])}
           \rho(f_z\circ g\circ f_y^{-1}).
\end{alignat*}
The fact that this satisfies $F_\rho(h)F_\rho(g)=\theta([h|g])
F_\rho(h\circ g)$ is a slightly brutal but straight forward
calculation, this becomes easier after the factor out-front has
been explained in Section~\ref{Subsubsection:Explanation} below.

\begin{lemma}
Suppose $\cG$ is a finite, connected groupoid, $\theta\in
Z^2(\cG,U(1))$ is a two-cycle, and $x$ is an object of $\cG$. Any
$\theta$--twisted representation $F\colon \pretheta\cG\to \Vect$
is determined by its restriction
$F|_x\colon{}^{\theta|_x}\!\Aut(x)\to \Aut(F(x))$ in the sense
that $F$ is isomorphic to the representation induced up from
$F|_x$ (for any choice of morphisms $\{f_y\colon y\to x\}_{y\ne
x}$).
\end{lemma}
\begin{proof}
We just need to construct a natural isomorphism $T\colon F\to
F_\rho$, for this we take $T(y):=F(f_y)$ then the fact that
$T(z)F(g)=F_\rho(g)T(y)$ follows when you observe that
$F(f_y)^{-1}=\theta([f_y|f_y^{-1}])^{-1}F(f_y^{-1})$ and
$\theta([f_y|f_y^{-1}])^{-1}\theta([f_z g|f_y^{-1}])
=\theta([f_zgf_y^-1|f_y])^{-1}$.
\end{proof}

The case of a non-connected groupoid is then an immediate
corollary:

\begin{thm} \label{Theorem:DecomposeRep}%
The twisted representation group of a finite groupoid
splits as the direct sum of the twisted representation groups of
its automorphism groups: in other words
$$\thetaRep(\cG)\cong \bigoplus_{[x]\subset\cG} \Rep^{\theta|x}(\!\Aut(x)),$$
where $x$ runs over one object from each connected component of
$\cG$.
\end{thm}


\subsubsection{Explanation}
\label{Subsubsection:Explanation} As explained in
Section~\ref{Section:Equivalence} the system of maps $\{f_y\colon
y\to x\}$ defines a ``retraction'' $r\colon \cG\to
\overline{\Aut(x)}$ which gives an inverse equivalence to the
inclusion $i\colon \overline{\Aut(x)}\to\cG$. The functor $r\circ
i$ is the identity, but the functor $i\circ r$ is only isomorphic
to the identity with a natural transformation $T\colon
\id_\cG\nattrans i\circ r$ being given by $T(y):=f_y$.

Given a genuine representation of $\rho\colon\Aut(x)\to \Aut(V)$
we can pull this back to $\cG$ to get $r^*\rho\colon \cG\to \Vect$
a representation of $\cG$.  If $\rho$ is the restriction
$F|_x=i^*F$ then the representation $r^*\circ i^*F$ is isomorphic
to the original representation via $T^*F$.

Unfortunately this argument does not work in the twisted case for
the following reason.  If $\psi$ is a two-cocycle on $\Aut(x)$ and
$\rho\colon {}^\psi\!\Aut(x)\to \Aut(V)$ is a $\psi$--twisted
representation then $r^*\rho$ is an $r^*\psi$--twisted
representation, so if $\psi=i^*\theta=\theta|_x$ then $r^*\rho$ is
an $(i\circ r)^*\theta$--twisted representation of $\cG$. However,
in general $(i\circ r)^*\theta$ is not \emph{equal} to $\theta$
but only \emph{cohomologous} to it, in other words they will
differ by the boundary of a one-cochain $\epsilon$, so $(i\circ
r)^*\theta=d\epsilon\cdot\theta$ (remember we write it
multiplicatively as the cochains take values in $U(1)$).
Fortunately there is a way around this. Using standard topological
methods we can write down a suitable $\epsilon$ by using the
natural transformation $T$. It is probably easiest to do it in
slightly more generality.
\begin{thm}
Suppose $\cG$ and $\cH$ are groupoids, $K,\check K\colon
\cG\to\cH$ are isomorphic functors, $T\colon K\nattrans\check K$
is any natural transformation, $\theta\in Z^2(\cH,U(1))$ is a
two-cocycle on $\cH$, and $\epsilon\in C^1(\cG,U(1))$ is the
one-cochain on $\cG$ defined by
  \[
    \epsilon(x_0\stackrel{g}\rightarrow x_1):=
       \frac{\theta([T(x_1)|K(g)])}
       {\theta([\check K(g)|T(x_0)])}.
  \]
Then the pull-backs of $\theta$ along $K$ and $\check K$ are
cohomologous via $\epsilon$, ie.\
  \[
    K^*(\theta)=d\epsilon\cdot\check{K}^*(\theta)\in
    Z^2(\cG,U(1)).
  \]
Furthermore, if $F\colon \pretheta\cH\to \Vect$ is a
$\theta$--twisted representation then $\epsilon\cdot
\check{K}^*(F)$ is a $K^*(\theta)$--twisted representation which
is isomorphic to $K^*(F)$ via the natural transformation $F\circ
T\colon K^*(F)\nattrans \epsilon\cdot \check{K}^*(F)$.
%
\end{thm}
\begin{proof}
  Consider the homotopy $H_T\colon \cB \cG\times I \to\cB\cH$
  given in Section~\ref{Section:Equivalence}.  In particular fix a
  two-simplex $\sigma$ in $\cB\cG$ and look at the prism
  $P:=H_T(\sigma\times  I )$ in $\cB \cH$.
 \begin{figure}[h]
   $$\vcenter{\relabelbox\small\epsfbox{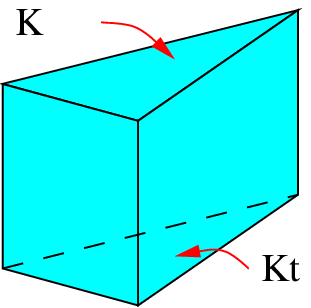}
 \relabel{K}{$\check K(\sigma)$}
 \relabel{Kt}{$K(\sigma)$}
 \endrelabelbox}$$
 \end{figure}

  \noindent Observe that the boundary of the prism is given by $\partial
  P=K(\sigma)-\check K(\sigma)+H_T(\partial \sigma \times  I )$.
  As $\theta$ is a two-cocycle, it is trivial on boundaries such as
  $\partial P$, thus
  \[
    \theta(K(\sigma))^{-1}\theta(\check K(\sigma))
    \theta(H_T(\partial \sigma \times  I ))=1,
  \]
  or in other words,
  \[K^*\theta(\sigma)=d\epsilon(\sigma)\cdot \check K^*\theta(\sigma)\]
  where $\epsilon$ is the one-cochain defined by
  $\epsilon (\hat\sigma):=\theta(H_T(\hat\sigma\times  I ))$
  for any one-chain $\hat \sigma$ on $\cG$ ---
  remember $d\epsilon(\sigma)$ is, by definition,
  $\epsilon(\partial\sigma)$.  Then we see immediately from the
  formula given in Section~\ref{Section:Equivalence} that for a
  morphism $g\colon x_0\to x_1$ that $H_T([g]\times
   I )=[T(x_1)|K(g)]-[\check K(g)|T(x_0)]$, from which the
  formula for $\epsilon$ follows.
\end{proof}

In the case that we are interested above, we have $\cH:=\cG$,
$K:=\id_\cG$, $\check K:=i\circ r$, so
 $\check K(y\stackrel g\to z):=f_z\circ g\circ f_y^{-1}$,
$T(y):=f_y$, from which we get $\check K^*(F)(y\stackrel g\to
z)=F|_x(f_z\circ g\circ f_y^{-1})$, and so we recover the formula
from the previous subsection.

\section{The Drinfeld double of a finite group}
We can now see how the general theory of the previous section
easily recovers many facts about the twisted Drinfeld double of a
finite group.  Indeed it is only through this point of view that I
understand the Drinfeld double.  Firstly we identify the twisted
Drinfeld double as the twisted groupoid algebra of the loop
groupoid of the finite group.  As the loop groupoid of a finite
group is the same as the conjugation action groupoid, in
Section~\ref{Subsection:TwistedKtheory} we immediately recover
that the representation group of the Drinfeld double is the
twisted equivariant K-theory of the group.  The fact that the
twisted representation group of a groupoid is the direct sum of
the representation groups of its automorphism groups translates in
Section~\ref{Subsection:DPRInduction} into the decomposition of
the representation group of the Drinfeld double \textit{\`a la}
Dijkgraaf-Pasquier-Roche induction.  Finally in
Section~\ref{Subsection:DrinfeldCharacterTheory} we see how the
characters of the Drinfeld double classify representations and are
`twisted elliptic characters' on the group.

\subsection{Identifying the twisted Drinfeld double}
For a finite group $G$ and a three-cocycle $\omega\in
Z^3_{\text{gp}}(G,U(1))$ Dijkgraaf, Pasquier and Roche \cite{DPR}
defined a quasi-Hopf algebra $D^\omega(G)$, the twisted Drinfeld
double of $G$.  Freed \cite{Freed:HigherAlgebraicStructures}
identified the untwisted Drinfeld double as the groupoid algebra
of the action groupoid $\cG_G(\Gconj)$.  I am ignoring the
coproduct, so we will just be considering $D^\omega(G)$ as an
associative algebra. Personally, I would take the following
theorem to be the definition of $D^\omega(G)$!
\begin{thm}
\label{Theorem:IdentifyDrinfeld}%
Given a finite group $G$ and a three-cocycle $\omega\in
Z^3_{\text{gp}}(G,U(1))$, let $\tau\colon
Z^3_{\text{gp}}(G,U(1))\to Z^2(\Lambda \overline G,U(1))$ be the
loop transgression map, then the twisted Drinfeld double
$D^\omega(G)$ of $G$ is isomorphic to the twisted groupoid algebra
${}^{\tau(\omega)}\!\C(\Lambda \overline G)$.
\end{thm}
\begin{proof}
This is just unpacking the definition of the twisted groupoid
algebra.  A basis of this algebra is indexed by the set of
morphisms in $\Lambda \overline G$, so we will write the basis
element correspond to $gxg^{-1}\stackrel g \leftarrow x$ for
$x,g\in G$ as $\langle \stackrel g \leftarrow x\rangle$. By
definition the product of two basis elements is given by
\begin{align*}
  \langle \stackrel h \leftarrow y\rangle \cdot
  \langle \stackrel g \leftarrow x\rangle
   &:=
   \delta_{y,gxg^{-1}}
       \tau(\omega)([\stackrel {h}\leftarrow gxg^{-1}|
                    \stackrel g \leftarrow x])~
       \langle \stackrel{hg}\longleftarrow x \rangle
       ,
\end{align*}
where
\[
\tau(\omega)([\stackrel{h}\leftarrow gxg^{-1}|
                   \stackrel g \leftarrow x]):=
\frac{\omega(h,g,x)\omega(hgx(hg)^{-1},h,g)}{\omega(h,gxg^{-1},g)}.
\]
What we have here ``up to some changes in notation''
\cite{Freed:HigherAlgebraicStructures} is the definition of
Dijkgraaf, Pasquier and Roche; actually the notation is
essentially the same as that of Freed
\cite{Freed:HigherAlgebraicStructures}.
\end{proof}
This identification has various consequences.

\subsection{Twisted K-theory} \label{Subsection:TwistedKtheory}%
The first consequence of the above theorem can be viewed as a baby
version of the Freed-Hopkins-Teleman theorem connecting the
Verlinde algebra to ${}^{\tau(\omega)}\!K_G(\Gconj)$ the twisted
equivariant K-theory of $G$.  The key here is to understand what
the definition of the twisted K-group is.

The untwisted equivariant K-group $K_G(\Gconj)$ of $G$ is defined
to be the Grothendieck group of equivariant vector bundles over
$\Gconj$, which is $G$ equipped with the conjugation action.
However, an equivariant bundle over $\Gconj$ is precisely the same
thing as a representation of the action groupoid $\cG_G(\Gconj)$;
hence $K_G(\Gconj)\cong \Rep(\cG_G(\Gconj))$. So by the twisted
equivariant K-theory ${}^{\theta}\!K_G(\Gconj)$ we mean the
twisted representation group $\thetaRep(\cG_G(\Gconj))$ for
$\theta$ some two-cocycle on the action groupoid $\cG_G(\Gconj)$.
Now by Section~\ref{Subsection:ActionGroupoid}, the latter
groupoid can be identified with the loop groupoid $\Lambda
\overline G$, so a two-cocycle on the conjugation action groupoid
is the same thing as a two-cocycle on the loop groupoid. In
particular if $\omega\in Z^3_{\text{gp}}(G,U(1))$ then
$\tau(\omega)\in Z^2(\Lambda\oG,U(1))$ so we can define  the
twisted equivariant K-theory ${}^{\tau(\omega)}\!K_G(\Gconj)$ to
be $\Rep^{\tau(\omega)}\!(\Lambda\oG)$.  We can then get the
following result.
\begin{prop}
The representation group of the twisted Drinfeld double of a
finite group $G$ is naturally isomorphic to the twisted
equivariant K-group of $G$:
  $$\Rep(D^\omega(G))\cong {}^{\tau(\omega)}\!K_G(\Gconj).$$
\end{prop}
\begin{proof}
  The twisted equivariant K-theory was defined to be
  $\Rep^{\tau(\omega)}\!(\Lambda\oG)$, but by
  Proposition~\ref{Prop:AlgGpdRep} this is the representation
  group of the twisted groupoid algebra
  ${}^{\tau(\omega)}\!\C\Lambda\oG$, and by
  Theorem\ref{Theorem:IdentifyDrinfeld} this algebra is just the
  twisted Drinfeld double $D^\omega(G)$.
\end{proof}
\subsection{Representation theory} \label{Subsection:DPRInduction}%
As the twisted Drinfeld double is a twisted groupoid algebra, the
representation theory decomposes according to
Theorem~\ref{Theorem:DecomposeRep} into the twisted representation
groups of the automorphism groups of the groupoid.  In the case of
the loop groupoid $\Lambda\oG$, an object is an element of $G$ and
the automorphism group of an object $x$ is the centralizer $C_x$.
Thus by Theorem~\ref{Theorem:IdentifyDrinfeld},
Theorem~\ref{Theorem:DecomposeRep}, and
Proposition~\ref{Prop:AlgGpdRep} we get following decomposition.
\begin{thm} The representation group of the twisted Drinfeld
double of a finite group $G$ decomposes as the twisted
representation groups of the centralizers of the conjugacy classes
of $G$:
\[\Rep(D^\omega(G))\cong \bigoplus_{[x]\subset G}
\Rep^{\tau(\omega)|_x}\!(C_x),
\]
where $x$ runs over one representative from each conjugacy class
of $G$.
\end{thm}
\noindent Given a twisted representation of the centralizer of an
element in $G$, the way to get a representation of the twisted
Drinfeld double is described in
Section~\ref{Subsubsection:DPRInduction} and this is immediately
seen to be precisely DPR induction \cite{DPR}.

\subsection{Character theory}
\label{Subsection:DrinfeldCharacterTheory}
%

The ``character'' of a representation $\rho\colon D^\omega(G)\to
\End(V)$ of the twisted Drinfeld double is the function
$\chi_\rho$ on pairs $(g,x)$ defined by
\[
  \chi_\rho(g,x):=\Tr\circ\rho
     \bigl(\langle \stackrel{g}{\leftarrow}x\rangle\bigr)
\]
and this will vanish unless $x=gxg^{-1}$ (or $xg=gx$ if you
prefer).  By Section~\ref{Section:TwistedCharacters} this is the
twisted character of the corresponding $\tau(\omega)$--twisted
representation of the loop groupoid $\Lambda\oG$, so it will be a
section of the line bundle $\tau^2(\omega)_\C$ over
$\Lambda^2\oG$.  Such a section is just a function $\chi$  on
pairs $(g,x)\in G\times G$ such that $xg=gx$ (so that $\chi(g,x)$
really means $\chi(x\stackrel g \leftarrow x)$) such that for any
$h\in G$
\begin{align*}
  \chi(hgh^{-1},hxh^{-1})
    &=\tau^2(\omega)\left([h](x\stackrel g \leftarrow x)\right) \chi(g,x)\\
    &=\frac{\tau(\omega)([hgh^{-1}|h]x)}{\tau(\omega)([h|g]x)} \chi(g,x)\\
    &=\frac{\omega(h,x,g)\omega(hgh^{-1},h,x)\omega(hxh^{-1},hgh^{-1},h)}
          {\omega(h,g,x)\omega(hxh^{-1},h,g)\omega(hgh^{-1},hxh^{-1},h)}
           \chi(g,x).
\end{align*}
We can call functions of this form \emph{$\omega$--twisted
elliptic characters of $G$}.  Thus we obtain a map
\[
  \chi\colon \Rep(D^\omega(G))\to
     \{\text{$\omega$--twisted elliptic characters of $G$}\}
\]
Theorem~\ref{Theorem:DimSectionsGerbe} gives that representations
of the Drinfeld double are classified by their character:
\begin{thm}
  The complexification of the character map is an isomorphism
\[
  \chi\colon \Rep(D^\omega(G))\otimes \C
     \stackrel{\cong}{\longrightarrow}
     \{\text{$\omega$--twisted elliptic characters of $G$}\}
\]
and an orthonormal basis of the $\omega$--twisted elliptic
characters of $G$ is given by the characters of the irreducible
representations of $D^\omega(G)$.
\end{thm}
Using Corollary~\ref{Cor:NumberTwistedReps} we can give the
dimension of the space of $\omega$--twisted elliptic characters,
which is thus the number of irreducible representations of the
twisted Drinfeld double.
\begin{thm}
  The number of irreducible representations of the twisted
  Drinfeld double $D^\omega(G)$ is given by
  \begin{align*}\rank \Rep(D^{\omega}(G))
    &=
    \int_{\Lambda^3\oG}\tau^3(\omega)
    = \frac{1}{|G|}\sum
     \frac{\omega(h,x,g)\omega(g,h,x)\omega(x,g,h)}
          {\omega(h,g,x)\omega(x,h,g)\omega(g,x,h)}
  \end{align*}
  where the sum is over triples $(h,g,x)\in G\times G\times G$
  which mutually commute.
\end{thm}

\begin{figure}
 $$\epsfbox{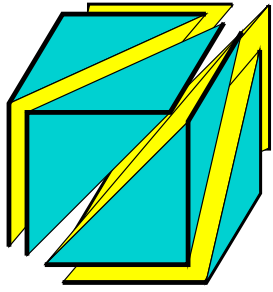}$$
  \caption{The three-torus decomposed.}
  \label{Figure:ThreeTorusCube}
\end{figure}

Actually this is immediately seen to be the invariant of the
three-torus, which is decomposed into six three-simplices as in
Figure~\ref{Figure:ThreeTorusCube}.

A simple example is given by taking $G=\Z_2\times\Z_2\times\Z_2$,
with generators $a$, $b$ and $c$ and with three-cocycle $\omega\in
Z^3(\Z_2\times\Z_2\times\Z_2,U(1))$ defined by
\[\omega(a^{u_a}b^{u_b}c^{u_c},a^{v_a}b^{v_b}c^{v_c},a^{w_a}b^{w_b}c^{w_c})
    =(-1)^{u_a.v_b.w_c}.\]
A quick calculation on {\tt maple} of the previous formula gives
  \[\rank \Rep(D^{\omega}(\Z_2\times\Z_2\times\Z_2))=22.\]
A little more investigation reveals that this is accounted for by
8 one-dimensional representations and 14 two-dimensional
representations.  A simple check on these numbers is that
$|\Z_2\times\Z_2\times\Z_2|^2=64=14\times 2^2 +8\times 1^2$.


\begin{thebibliography}{20}
\bibitem{AltschulerCoste:QuasiQuantum}
   D.~Altschuler and A.~Coste,
   \textit{Quasi-quantum groups, knots, three-manifolds, and topological field theory},
   Comm.\ Math.\ Phys.\ \textbf{150} (1992), no.~1,
   83--107.  \texttt{arXiv:hep-th/9202047}
\bibitem{AltschulerCosteMaillard}
   D.~Altschuler, A.~Coste and J-M.~Maillard,
  \textit{Representation theory of twisted group double},
  \texttt{arXiv:hep-th/030925}.
\bibitem{Benson:BookII} D.~J.~Benson,
  Representations and cohomology. II. Cohomology of groups and modules,
  Cambridge Studies in Advanced Mathematics, \textbf{31},
  Cambridge University Press, 1991.
\bibitem{BrylinskiMcCaughlin} J.-L.~Brylinski and D.~McLaughlin,
    \textit{The geometry of degree-four characteristic classes
     and of line bundles on loop spaces. I},
    Duke Math.\ J.\ \textbf{75} (1994), no.~3, 603--638.
\bibitem{DPR} R~Dijkgraaf, V~Pasquier and P.~Roche,
   \textit{Quasi Hopf algebras, group cohomology and orbifold models},
   Nuclear Physics B (Proc.~Suppl) \textbf{18B} (1990) 60--72.
\bibitem{DijkgraafWitten} R.~Dijkgraaf, E.~Witten,
  \textit{Topological gauge theories and group cohomology},
  Comm.\ Math.\ Phys.\ \textbf{129} (1990), no. 2, 393--429.
\bibitem{Freed:HigherAlgebraicStructures} D.~Freed,
   \textit{Higher algebraic structures and quantization},
   Comm.\ Math.\ Phys.\ \textbf{159} (1994), no.~2, 343--398.
   \texttt{arXiv:hep-th/9212115}
\bibitem{FreedQuinn} D.~Freed and F.~Quinn,
  \textit{Chern-Simons theory with finite gauge group},
  Comm.\ Math.\ Phys.\ \textbf{156} (1993), no.\ 3, 435--472.
  \texttt{arXiv:hep-th/9111004}
\bibitem{Gajer:DeligneCohomology} P.~Gajer,
  \textit{Geometry of Deligne cohomology},
  Invent.\ Math.\ \textbf{127} (1997), no.\ 1, 155--207.
  \texttt{arXiv:alg-geom/9601025}
\bibitem{GomiTeramshima:FiberIntegration} K.~Gomi and
   Y.~Terashima,
   \textit{A fiber integration formula for the smooth Deligne
   cohomology},
   Int.\ Math.\ Res.\ Notices \textbf{13} (2000) 699--708.
\bibitem{Karpilovsky:Book} G.~Karpilovsky,
  Projective representations of finite groups,
  Monographs and Textbooks in Pure and Applied Mathematics,
  \textbf{94}, Marcel Dekker, Inc., New York, 1985.
\bibitem{LupercioUribe:GerbesTwistedKtheory}
  E.~Lupercio and B.~Uribe,
  \textit{Gerbes over orbifolds and twisted $K$-theory},
  Comm.\ Math.\ Phys.\ \textbf{245} (2004), no.~3, 449--489.
  \texttt{arXiv:math.AT/0105039}
\bibitem{Strickland} N.~P.~Strickland,
  \textit{$K(N)$-local duality for finite groups and groupoids},
  Topology \textbf{39} (2000), no.~4, 733--772.
  \texttt{arXiv:math.AT/0011109}
\end{thebibliography}
\end{document}